\documentclass[]{interact}

\usepackage{epstopdf}
\usepackage[caption=false]{subfig}

\usepackage[numbers,sort&compress]{natbib}
\bibpunct[, ]{[}{]}{,}{n}{,}{,}

\theoremstyle{plain}
\newtheorem{theorem}{Theorem}[section]
\newtheorem{lemma}[theorem]{Lemma}
\newtheorem{corollary}[theorem]{Corollary}
\newtheorem{proposition}[theorem]{Proposition}

\theoremstyle{definition}
\newtheorem{definition}[theorem]{Definition}
\newtheorem{example}[theorem]{Example}

\theoremstyle{remark}
\newtheorem{remark}{Remark}

\usepackage{amsmath, amssymb}
\usepackage{mathtools}
\usepackage{enumitem}
\usepackage{hyperref}
\usepackage{ulem}
\usepackage{float}   
\floatstyle{plain}
\usepackage{graphicx}

\DeclareMathOperator{\conv}{conv}
\DeclareMathOperator{\cl}{cl}

\DeclareMathOperator*{\argmin}{argmin}

\newcommand{\norm}[1]{\left\lVert#1\right\rVert}
\newcommand\numberthis{\addtocounter{equation}{1}\tag{\theequation}}
\usepackage{changes}

\usepackage{tikz}
\usetikzlibrary{matrix,arrows,decorations, arrows.meta, decorations.markings}
\usepackage{pgfplots}
\usepgfplotslibrary{fillbetween}

\newenvironment{MSCcodes}{%
        \par\addvspace{13pt plus2pt minus1pt}
  \keywordfont\noindent{\bfseries MATHEMATICS SUBJECT CLASSIFICATION (2020)\\}\ignorespaces%
}{%
        \par\addvspace{13pt plus2pt minus1pt}
  \@endparenv}
  
\usepackage{fancyhdr}
\begin{document}

\title{A new problem qualification based on approximate KKT conditions for Lipschitzian optimization with application to bilevel~programming}

\author{
\name{Isabella K{\"a}ming\textsuperscript{a}\thanks{Contact: Isabella K{\"a}ming. Email: isabella.kaeming@tu-dresden.de}, Andreas Fischer\textsuperscript{a} and Alain B. Zemkoho\textsuperscript{b}}
\affil{\textsuperscript{a} Institute of Numerical Mathematics, Technische Universit\"{a}t Dresden, 01062 Dresden, Germany\\ \textsuperscript{b} School of Mathematics, University of Southampton, Southampton SO17 1BJ, UK}
}

\maketitle

\begin{abstract}
When dealing with general Lipschitzian optimization problems, there are many problem classes where even weak constraint qualifications fail at local minimizers. In contrast to a constraint qualification, a problem qualification does not only rely on the constraints but also on the objective function to guarantee that a local minimizer is a Karush-Kuhn-Tucker (KKT) point. For example, calmness in the sense of Clarke is a problem qualification. In this article, we introduce the Subset Mangasarian-Fromovitz Condition (subMFC). This new problem qualification is derived by means of a nonsmooth version of the approximate KKT conditions, which hold at every local minimizer without further assumptions. A comparison with existing constraint and problem qualifications reveals that subMFC is strictly weaker than quasinormality and can hold even if the local error bound condition, the cone-continuity property, the Guignard constraint qualification and calmness are violated. Furthermore, we emphasize the power of the new problem qualification within the context of bilevel optimization. More precisely, under mild assumptions on the problem data, we suggest a version of subMFC that is tailored to the lower-level value function reformulation. It turns out that this new condition can be satisfied even if the widely used partial calmness condition does not hold.
\end{abstract}



\section{Introduction}
We consider mathematical optimization problems of the form
\begin{equation}
\min_{x} f(x) \quad\text{s.t.}\quad x \in D\coloneqq \{ x\in\mathbb{R}^n \mid g(x)\leq 0\} \tag{P}\label{P} 
\end{equation}
with locally Lipschitz continuous functions $f:\mathbb{R}^n \to \mathbb{R}$ and $g:\mathbb{R}^n\to \mathbb{R}^q$.
For simplicity, we only regard problems with inequality constraints in the course of this paper. Later on, we will see by examples that equality constraints can be handled as well by splitting them into two inequality constraints. It is known that a local minimizer of problem \eqref{P} satisfies the KKT conditions if some suitable constraint qualification (CQ) is fulfilled.
Unfortunately, there are situations where even weak CQs fail. In such cases we are faced with the problem that it is hard to analyze whether we can expect the KKT conditions to hold at local minimizers of the problem. For example, this difficulty generically occurs in bilevel optimization \cite{JY3} or for mathematical programs with complementarity constraints (MPCCs) \cite{RA7,AR1}. 

Therefore, instead of directly considering the KKT conditions, we are interested in optimality conditions that do not require the fulfillment of a CQ. For the smooth counterpart of problem~\eqref{P}, an established concept is the usage of sequential optimality conditions, as they are valid independently of CQs, see \cite{RA3,RA4,RA7,EB1}.
Moreover, sequential optimality conditions imply the KKT conditions under specially tailored conditions which are weak enough such that they can actually be satisfied \cite{RA6,RA5}. Thus, the procedure of combining sequential optimality conditions with a weak condition is also of interest to obtain optimality conditions for instances of the aforementioned problem classes, and has been carried out in \cite{RA7,AR1} for MPCCs with continuously differentiable problem data. 
Another advantage of sequential optimality conditions is that they can be used to justify termination criteria for numerical algorithms, as done for augmented Lagrangian methods \cite{EB1}, SQP methods \cite{LQ1}, inexact restoration methods \cite{LB1,AF1} and others, like \cite{RA3,RA4,JM1}. 

In this paper, as sequential optimality conditions, we consider the well-known approximate KKT (AKKT) conditions and transfer the definition used in \cite{LQ1} to our nonsmooth setting. 
The main focus is then to introduce and analyze a new condition which guarantees that a point fulfilling the AKKT conditions is already a KKT point. 
One could therefore be tempted to assume that, in terms of \cite{EB1}, this new condition belongs to the class of \textit{strict} CQs. 
Strict CQs are all CQs which ensure that an AKKT point is also a KKT point. 
For example, the class of strict CQs contains CQs like the linear independence CQ (LICQ), Mangasarian-Fromovitz CQ (MFCQ), constant rank CQ (CRCQ), constant positive linear dependence condition (CPLD) and the cone-continuity property (CCP) \cite{RA5}. 
The latter is the weakest strict CQ \cite{RA5,EB1} and also called asymptotic regularity in a nonsmooth context \cite{PM1}. 
In contrast to this, the new condition we are interested in does \textit{not} belong to the class of strict CQs. 
In fact, it directly depends on the specific sequence fulfilling the AKKT conditions for the given point, which is a significant difference to strict CQs that are stated completely independent of AKKT conditions. 
Due to the latter circumstance, the new condition is \textit{not even} a CQ, as the AKKT conditions and therefore the condition itself involve the objective function. 
Hence, we refer to such conditions as \textit{problem qualifications}. Similarly to CQs, problem qualifications also guarantee that the KKT conditions hold at a local minimizer. A well-known problem qualification is calmness in the sense of Clarke \cite[Definition 6.4.1]{FC1}.

It turns out that the proposed novel coupling of the AKKT conditions with a directly corresponding condition provides a very powerful problem qualification for \eqref{P}, which we call \textit{Subset Mangasarian-Fromovitz Condition (subMFC)}. 
Indeed, exploiting the information obtained by the specific fulfilled AKKT conditions enables subMFC to still maintain the most important property of strict CQs, which is serving as a condition to ensure that an AKKT point is already a KKT point. 
At the same time, it can be shown that subMFC can hold even if CCP (the weakest strict CQ) fails. 
Besides this advantage, subMFC has the following additional benefits:
\begin{itemize}
\item It is easy to understand, since it can be interpreted as a weakened version of the MFCQ.
\item A comparison with existing CQs for \eqref{P} reveals that subMFC is a quite weak condition that can be satisfied by a wide range of optimization problems. In particular, subMFC is strictly implied by quasinormality. Additionally, examples show that subMFC is independent of the local error bound condition and the Guignard CQ, and that subMFC can hold even if the problem qualification calmness in the sense of Clarke is violated. 
\newpage
\item It can be handled easier than CCP, which is hard to verify \cite{RA5} since it depends on \textit{all} subgradients that occur in the neighborhood of the considered point. With subMFC, generally significantly fewer subgradients need to be considered.
\item SubMFC can be used to treat all problem classes covered by \eqref{P}, which includes MPCCs \cite{ZL1}, mathematical programs with vanishing constraints \cite{WA1} or mathematical programs with switching constraints \cite{PM4}.
\end{itemize}
As for the last advantage, in this paper, we specifically demonstrate the power of our approach by applying it to bilevel optimization problems. For this purpose, let us consider the standard optimistic bilevel optimization problem 
\vspace{-0.1em}
\begin{equation}
\min_{x,y} F(x,y) \quad\text{s.t.}\quad G(x,y)\leq 0,\quad y\in S(x) \label{BP}\tag{BP}
\vspace{-0.2em}
\end{equation}
with the lower-level solution set
\[ 
S(x)\coloneqq \argmin_y \{ \hat{f}(x,y) \mid \hat{g}(x,y)\leq 0 \}, 
\vspace{-0.5em}
\]
where the functions $F:\mathbb{R}^n \times \mathbb{R}^m \to \mathbb{R}$, $G:\mathbb{R}^n \times \mathbb{R}^m\to \mathbb{R}^p$, $\hat{f}:\mathbb{R}^n \times \mathbb{R}^m\to \mathbb{R}$ and $\hat{g}:\mathbb{R}^n \times \mathbb{R}^m\to \mathbb{R}^r$ are assumed to be locally Lipschitz continuous. 
A widely used single-level reformulation of \eqref{BP} is based on the lower-level value function $\varphi:\mathbb{R}^n \to \mathbb{R}$, which is defined by
\vspace{-0.7em}
\[ 
\varphi(x) \coloneqq \min_y \{ \hat{f}(x,y) \mid \hat{g}(x,y) \leq 0 \}. 
\vspace{-0.2em}
\]
To ensure that $\varphi$ is well-defined, we assume throughout that the lower-level solution set $S(x)$ is not empty for each $x\in\mathbb{R}^n$. 
The lower-level value function reformulation is then given as
\vspace{-0.3em}
\begin{equation}
\min_{x,y} F(x,y) \quad\text{s.t.}\quad G(x,y)\leq 0,\quad \hat{g}(x,y)\leq 0,\quad \hat{f}(x,y) - \varphi(x) \leq 0. \label{LLVF}\tag{LLVFR} 
\vspace{-0.2em}
\end{equation}
By definition of $\varphi$, this problem is equivalent to \eqref{BP} concerning global and local solutions.
If $\varphi$ is locally Lipschitz continuous, it is clear that \eqref{LLVF} is a special case of \eqref{P}.~Thus, we will apply our results related to subMFC to derive the new problem qualification \textit{LLVFR-subMFC} that is tailored to \eqref{LLVF}. The major advantages of LLVFR-subMFC are:
\begin{itemize}
\item LLVFR-subMFC can hold even if partial calmness is violated. Partial calmness is a condition which has emerged as a standard condition for deriving KKT-type optimality conditions for bilevel optimization problems treated by the lower-level value function reformulation \cite{SD2,SD1,SD3,AF2,JY3}.
\item In contrast to most of the existing literature, we do not require any additional assumptions on the bilevel problem besides the local Lipschitz continuity of $F$, $G$, $\hat{f}$ and $\hat{g}$ and the existence and local Lipschitz continuity of $\varphi$. In particular, the given functions do not necessarily need to be differentiable or convex in any variable and coupling constraints in the upper level are allowed.
\end{itemize}
The remainder of the paper is structured as follows. In the next section, we recall some required basic terminology.
In Section 3, we discuss the concept of approximate KKT points in a Lipschitzian setting which is the basis for the introduction of the new problem qualification subMFC. The main result in this section shows that any feasible point of \eqref{P} is a KKT point provided that it satisfies subMFC. In Section 4, we focus on the comparison of subMFC to existing CQs and the problem qualification calmness. Afterwards, the application to bilevel optimization problems is presented in Section~5. Final remarks are given in the last section of this paper.

\section{Preliminaries}
Throughout, we use $\norm{.}$ to indicate the Euclidean norm. The closed ball centered at $x\in\mathbb{R}^n$ with radius $\epsilon>0$ is given by $B(x,\epsilon) \coloneqq \{ y\in\mathbb{R}^n \mid \norm{y-x} \leq \epsilon\}$.  For \mbox{$x=(x_1,\dotsc,x_n)^\top \in\mathbb{R}^n$}, we define $x_+ \coloneqq (\max\{0,x_1\},\dotsc,\max\{0,x_n\})^\top$, whereas $d_C(x)$ stands for the Euclidean distance of $x$ to a set $C\subseteq\mathbb{R}^n$. To shorten notation, from now on a tuple $(x_1,\dotsc,x_n)$ of real numbers $x_1,\dotsc,x_n\in\mathbb{R}$ will be identified with the vector $(x_1,\dotsc,x_n)^\top \in\mathbb{R}^n$.
Whenever a function $\phi:\mathbb{R}^n\to\mathbb{R}$ is differentiable at $\bar{x}\in\mathbb{R}^n$, we denote the gradient of $\phi$ at $\bar{x}$ by $\nabla\phi(\bar{x})$. For a function $\phi: \mathbb{R}^n \to \mathbb{R}$ that is Lipschitz continuous near $\bar{x}\in\mathbb{R}^n$, we use the well-known generalization of the gradient by Clarke \cite{FC1}. More in detail, with the \textit{generalized directional derivative} of $\phi$ at $\bar{x}$ in direction $d\in\mathbb{R}^n$ given by
\[ 
\phi^\circ(\bar{x};d) \coloneqq \limsup_{x\to \bar{x}, h\downarrow 0} \frac{\phi(x+hd)-\phi(x)}{h}, 
\]
the \textit{Clarke subdifferential} of $\phi$ at $\bar{x}$ is the set 
\[ 
\partial\phi(\bar{x}) \coloneqq \left\{ s\in\mathbb{R}^n \bigm| s^\top d \leq \phi^\circ(\bar{x};d)\,\,\, \forall d\in\mathbb{R}^n\right\}. 
\]
If $\phi$ is continuously differentiable, the Clarke subdifferential reduces to the gradient, i.e., $\partial\phi(\bar{x})=\{ \nabla \phi(\bar{x})\}$.
The following proposition provides useful properties of the Clarke subdifferential \cite[Proposition 2.1.2, Proposition 2.1.5]{FC1}.

\begin{proposition}
\label{prop:eig_clarke_subdiff}
Let $\phi:\mathbb{R}^n\to \mathbb{R}$ be locally Lipschitz continuous and $\bar{x}\in\mathbb{R}^n$. 
\begin{enumerate}[label=(\alph*)]
\item \label{propa} The set $\partial\phi(\bar{x})$ is nonempty, convex and compact.
\item \label{propb} If the sequence $\{ x^k \}_{k=1}^\infty$ converges to $\bar{x}$ and $s^k\in\partial \phi(x^k)$ for all $k\in\mathbb{N}$, then all accumulation points of the sequence $\{ s^k \}_{k=1}^\infty$ belong to $\partial \phi(\bar{x})$.
\item \label{propc} The map $x \mapsto \partial \phi(x)$ is upper semicontinuous at $\bar{x}$, i.e., for all $\epsilon>0$ there exists $\delta>0$ such that 
$
\partial \phi(x) \subseteq \partial \phi(\bar{x}) + B(x,\epsilon)$ holds for all $x\in B(\bar{x},\delta). 
$
\end{enumerate} 
\end{proposition} 

Other generalizations of the gradient include the regular and the limiting subdifferential \mbox{\cite{BM1,RR1}}. In detail, let $\phi:\mathbb{R}^n\to\mathbb{R}\cup\{-\infty,\infty\}$ be a lower semicontinuous function that is finite at $\bar{x}\in\mathbb{R}^n$. Then, the \textit{regular/Fr\'{e}chet subdifferential} of $\phi$ at $\bar{x}$ is the set 
\[ \widehat{\partial} \phi(\bar{x}) \coloneqq \left\{ s\in\mathbb{R}^n \bigm| \liminf_{x\to\bar{x},\,x\neq\bar{x}} \frac{\phi(x)-\phi(\bar{x})-\langle s,x-\bar{x} \rangle}{\norm{x-\bar{x}}} \geq 0 \right\} \]
and the \textit{limiting/Mordukhovich/basic subdifferential} of $\phi$ at $\bar{x}$ is the set 
\begin{align*}
\partial^L \phi(\bar{x}) \coloneqq \{ \bar{s} \in\mathbb{R}^n \mid \exists\, \{ (x^k, s^k) \}_{k=1}^\infty \subset \mathbb{R}^n \times \mathbb{R}^m:\, &(x^k,s^k) \to (\bar{x},\bar{s}), \phi(x^k)\to\phi(\bar{x}),\\ &s^k\in \widehat{\partial} \phi(x^k)\,\,\, \forall k\in\mathbb{N}\}.
\end{align*} 
Note that the inclusion $\widehat{\partial}\phi(\bar{x}) \subseteq \partial^L \phi(\bar{x}) \subseteq \partial\phi(\bar{x})$ holds.
For a point $\bar{x}\in\mathbb{R}^n$ that is feasible for \eqref{P}, we define the index set of active constraints as 
\[ 
I(\bar{x}) \coloneqq \{ j\in \{1, \dots, q\} \mid g_j(\bar{x}) = 0 \}. 
\] 
The nonsmooth extension of the classical MFCQ demands the positive linear independence of the subgradients of all constraints that are active in $\bar{x}$. Thus, its formal definition can be stated as follows.

\begin{definition}
\label{def:mfcq}
\textit{MFCQ} is said to be satisfied at a point $\bar{x}$ feasible for \eqref{P} if $I(\bar{x})=\emptyset$ or if, for all $v\geq 0$, $v\neq 0$, it holds that
\[ 
0 \displaystyle \notin \sum_{j\in I(\bar{x})} v_j \partial g_j(\bar{x}).
\]
\end{definition}

\section{The Subset Mangasarian-Fromovitz Condition}
\label{sec:subMFC}
In this section, we first introduce a suitable sequential optimality condition for the nonsmooth problem \eqref{P}, which, later on, will constitute a component of subMFC. To this end, let us recall the definitions of Fritz John (FJ) points and KKT points \cite[Theorem 6.1.1]{FC1}. 

\begin{definition}
\label{def:fj}
A point $\bar{x}\in\mathbb{R}^n$ is called \textit{FJ point} of \eqref{P} if there exist multipliers $a\geq 0$, $u\in\mathbb{R}^{q}$, $(a,u)\neq (0,0)$ such that 
\vspace{-1em}
\begin{gather*}
 0 \in a \partial f(\bar{x}) + \sum_{j=1}^q u_j \partial g_j(\bar{x}),  \\ 
 u \geq 0,\quad g(\bar{x}) \leq 0,\quad u^\top g(\bar{x}) = 0.
\end{gather*} 
\end{definition}

\begin{definition}
\label{def:kkt}
A point $\bar{x}\in\mathbb{R}^n$ is called \textit{KKT point} of \eqref{P} if there exists a multiplier $u\in\mathbb{R}^{q}$ such that 
\vspace{-1em}
\begin{gather*}
 0 \in \partial f(\bar{x}) + \sum_{j=1}^q u_j \partial g_j(\bar{x}), \\ 
 u \geq 0,\quad g(\bar{x}) \leq 0,\quad u^\top g(\bar{x}) = 0.
\end{gather*} 
\end{definition}

It is known that a local minimizer $\bar{x}$ of problem \eqref{P} fulfills the KKT conditions if some CQ is satisfied at $\bar{x}$. But, on the one hand, CQs are often violated, making it challenging to figure out whether the KKT conditions hold at $\bar{x}$. Hence, we instead want to utilize optimality conditions that hold without further requirements. One could therefore be tempted to work with the FJ conditions, as they hold without further assumptions. But, in this case, the multiplier $a$ belonging to the objective function is allowed to vanish, such that the function being minimized is possibly not involved. On the other hand, even if we can expect the KKT conditions to hold, it is generally not possible to compute a KKT point of problem \eqref{P} exactly. Hence, numerical algorithms for the treatment of \eqref{P} typically incorporate termination criteria that test whether the current iterate is somehow close to a KKT point. Specifically, it is evaluated whether the KKT conditions given in Definition \ref{def:kkt} are approximately satisfied at the current iterate. This leads to the definition of approximate KKT points, which are well-suited to overcome all the aforementioned drawbacks.

\begin{definition}
\label{def:akkt}
A point $\bar{x}\in\mathbb{R}^n$ is called \textit{approximate KKT (AKKT) point} of \eqref{P} if there exists a 
sequence $\{(x^k,u^k, \delta^k, \epsilon^k)\}_{k=1}^\infty$ $\subset \mathbb{R}^n \times \mathbb{R}^q \times \mathbb{R}^q \times \mathbb{R}^n$ such that we have
\vspace{-0.4em}
\begin{gather}
 \epsilon^k \in \partial f(x^k) + \sum_{j=1}^q  u^k_j \partial g_j(x^k), \label{eq:akkt1}\\
 u^k \geq 0,\quad g(x^k) \leq \delta^k,\quad (u^k)^\top (g(x^k)-\delta^k) = 0 \label{eq:akkt2} 
\end{gather} 
for all $k\in\mathbb{N}$ as well as the convergence
\begin{equation}
 (x^k, \delta^k,\epsilon^k) \to (\bar{x},0,0) \quad \text{as}\,\,\, k\to\infty. \label{eq:akkt5} 
 \end{equation}
\end{definition}

Obviously, every KKT point is also an AKKT point. Definition \ref{def:akkt} can be seen as the nonsmooth counterpart of the definition in \cite{LQ1}. Note that different definitions of an AKKT point are used in the literature. For example, it can be verified that the definition above is equivalent to the one used in \cite{PM1}. If we additionally require $\delta^k\geq 0$, the definition obtained in this case would be equivalent to the one used in \cite{EB1} for the smooth case.

Next, we are going to prove that being an AKKT point is a necessary optimality condition that holds without further assumptions. Related results for nonsmooth settings can be found in \cite{PM2,PM1}, where the limiting subdifferential $\partial^L(\cdot)$ is used as a generalization of the gradient for nonsmooth functions. To maintain the paper self-contained, we now state an explicit proof that is tailored to our definition of an AKKT point and that is based on the Clarke subdifferential. Similar proof techniques for related results can be found in \cite{EB1,PM1}.
\begin{theorem}
\label{thm:locmin_akkt}
If $\bar{x}$ is a local minimizer of \eqref{P}, then $\bar{x}$ is an AKKT point of \eqref{P}.
\end{theorem}

\begin{proof}
Recall that $D$ denotes the feasible set of problem \eqref{P}. As $f$ and $g$ are locally Lipschitz continuous on $\mathbb{R}^n$ and $\bar{x}$ is a local minimizer of \eqref{P}, there exists $\epsilon>0$ such that $f$ and $g$ are locally Lipschitz continuous on $B(\bar{x},\epsilon)$ and $f(x)\geq f(\bar{x})$ holds for all $x\in D\cap B(\bar{x},\epsilon)$.
For $k\in\mathbb{N}$, let us consider the problem
\begin{equation}
\numberthis\label{eq:locmin_prob}
\min_x f(x) + k \norm{ g(x)_+ }^2 + \norm{ x-\bar{x} }^2\quad\text{s.t.}\quad x\in B(\bar{x},\epsilon). 
\end{equation}
Because a locally Lipschitz continuous objective function is to be minimized on a compact domain, the problem has a solution $x^k \in B(\bar{x},\epsilon)$ for all $k\in \mathbb{N}$ and the sequence $\{ x^k \}_{k=1}^\infty$ stays bounded. Therefore, without loss of generality, we may assume that $x^k \to x^\star$ holds for some $x^\star \in B(\bar{x},\epsilon)$ as $k \to \infty$.
By the definition of $x^k$ and the feasibility of $\bar{x}$, it follows that
\begin{equation}
\numberthis\label{eq:locmin}
f(x^k) + k \lVert g(x^k)_+ \rVert^2 + \lVert x^k-\bar{x} \rVert^2 \leq f(\bar{x}) + k \norm{ g(\bar{x})_+ }^2 + \norm{ \bar{x}-\bar{x} }^2 = f(\bar{x}) 
\end{equation}
is satisfied for all $k\in\mathbb{N}$. As $f$ is locally Lipschitz continuous on $B(\bar{x},\epsilon)$, $\{ f(x^k) \}_{k=1}^\infty$ is a bounded sequence. Thus, dividing by $k$ and passing to the limit in (\ref{eq:locmin}) yields $\norm{ g(x^\star)_+ }^2 = 0$, which is equivalent to $x^\star \in D$. From \eqref{eq:locmin}, it follows in particular that
\[ 
f(x^k) + \lVert x^k-\bar{x} \rVert^2 \leq f(\bar{x}) 
\]
holds for all $k\in\mathbb{N}$.
Using the latter, we obtain
\[
f(x^\star) \leq f(x^\star) + \norm{ x^\star-\bar{x} }^2 
= \lim_{k\to\infty} \left( f(x^k) + \lVert x^k-\bar{x} \rVert^2 \right) 
\leq f(\bar{x})
\leq f(x^\star),
\]
and hence, $x^\star = \bar{x}$. Thus, $x^k \to\bar{x}$ as $k \to \infty$ implying that $x^k$ lies in the interior of $B(\bar{x},\epsilon)$ for all $k\in\mathbb{N}$ large enough. Since $x^k$ is a solution to (\ref{eq:locmin_prob}), we obtain, using $f_1(x) \coloneqq \norm{x-\bar{x} }^2$ for all $x\in\mathbb{R}^n$, by \cite[Proposition 2.3.2]{FC1} that 
\[ 
0 \in \partial( f + k \norm{ g_+ }^2 + f_1 )(x^k).
\]
With the sum rule \cite[Proposition 2.3.3]{FC1} and the chain rule \cite[Theorem 2.3.9]{FC1} for the Clarke subgradient, we get
\begin{align*}
2(\bar{x}-x^k) \in \partial f(x^k) + \sum_{j:g_j(x^k)>0} 2k g_j(x^k) \partial g_j(x^k).
\end{align*} 
By setting
\[ 
\epsilon^k \coloneqq 2(\bar{x}-x^k), u^k_j \coloneqq \begin{cases} 2k g_j(x^k), & \text{if } g_j(x^k) > 0, \\ 0, & \text{if } g_j(x^k) \leq 0 \end{cases} \,\,\text{and}\,\, \delta^k_j \coloneqq \begin{cases} g_j(x^k), & \text{if }g_j(x^k) > 0, \\ \frac{1}{k}, & \text{if } g_j(x^k) \leq 0, \end{cases}  
\]
it can be checked that the conditions (\ref{eq:akkt1})--(\ref{eq:akkt5}) hold for the sequence $\{(x^k,u^k,\delta^k, \epsilon^k)\}_{k=1}^\infty$, and hence, $\bar{x}$ is an AKKT point.
\end{proof}

Thus, as it is possible to show that every local minimizer of problem \eqref{P} is the limit of a sequence of points which satisfy the AKKT conditions when the tolerances approach zero, the procedure of utilizing AKKT conditions for termination criteria in numerical algorithms is theoretically justified.

\begin{remark}
\label{rem:delta_greater_0}
The proof of Theorem \ref{thm:locmin_akkt} shows in particular that for a local minimizer $\bar{x}$ of \eqref{P}, it is always possible to find a sequence  $\{(x^k,u^k, \delta^k, \epsilon^k)\}_{k=1}^\infty$ with $\delta^k_j>0$ for all $j=1,\dotsc,q$ and $k\in\mathbb{N}$ for which Definition \ref{def:akkt} is fulfilled.
\end{remark}

Next, we want to proceed with the introduction of subMFC. To this end, for any $x\in\mathbb{R}^n$ that fulfills $g(x) \leq \delta$, we define the index set
\[ 
I(x,\delta) \coloneqq \{ j\in \{1, \dots, q\} \mid g_j(x) = \delta_j \}. 
\] 

\begin{remark}
\label{rem:seq_const}
Whenever an AKKT point with a corresponding sequence $\{(x^k,u^k, \delta^k, \epsilon^k)\}_{k=1}^\infty$ is considered in the following, we may assume without loss of generality that the associated set of active constraints $I(x^k,\delta^k)$ is constant for all $k\in\mathbb{N}$, which is reasonable as it is always possible to consider a suitable infinite subsequence. 
\end{remark}

\begin{definition}
\label{def:subMFC}
Let $\bar{x}$ be feasible for \eqref{P}. We say that the \textit{Subset Mangasarian-Fromovitz Condition (subMFC)} holds at $\bar{x}$ if there exists $I \subseteq I(\bar{x}) $ such that the following two conditions are satisfied:
\begin{enumerate}[label=(\roman*)]
\item \label{subMFCI} If $I\neq\emptyset$, it holds for all $v \geq 0, v\neq 0$ that
\[ 
0 \displaystyle \notin \sum_{j\in {I}} v_j \partial g_j(\bar{x}). 
\]
\item \label{subMFCII} There exists a sequence $\{(x^k,u^k,\delta^k,\epsilon^k)\}_{k=1}^\infty$ that fulfills system \eqref{eq:akkt1}--\eqref{eq:akkt5} with
\[ 
{I} = I(x^k,\delta^k). 
\]
\end{enumerate}
\end{definition}

By Theorem \ref{thm:locmin_akkt}, for the special case of $\bar{x}$ being a local minimizer of \eqref{P} there always exists at least one sequence that is a suitable candidate for part \ref{subMFCII} in the definition above, as part \ref{subMFCII} means that the considered point $\bar{x}$ is an AKKT point. Because of $I \subseteq I(\bar{x})$, part \ref{subMFCI} of the definition can be interpreted as a weakened MFCQ at $\bar{x}$, since we demand the positive linear independence of subgradients of active constraints only for some \textit{subset} of all constraints that are active in $\bar{x}$. However, the positive linear independence of these fewer subgradients already suffices to prove that the AKKT point satisfying \ref{subMFCII} is already a KKT point, as the following theorem shows.

\begin{theorem}
\label{thm:bqc_implies_kkt}
If $\bar{x}$ is feasible for \eqref{P} and subMFC holds at $\bar{x}$, then $\bar{x}$ is a KKT point of~\eqref{P}.
\end{theorem}

\begin{proof}
By Definitions \ref{def:akkt} and \ref{def:subMFC}, it follows that there is a sequence $\{(x^k,u^k,\delta^k,\epsilon^k)\}_{k=1}^\infty$ such that \eqref{eq:akkt2}, $(x^k, \delta^k,\epsilon^k) \to (\bar{x},0,0)$ for $k\to\infty$ and
\begin{equation}
\epsilon^k \in \partial f(x^k) + \sum_{j\in {I}}  u^k_j \partial g_j(x^k) \label{eq:1.7} \numberthis 
\end{equation}
hold with $I=I(x^k,\delta^k)$, taking into account Remark~\ref{rem:seq_const}. If $I=\emptyset$, Proposition \ref{prop:eig_clarke_subdiff}\ref{propb} and the local Lipschitz continuity of $f$ imply $0\in \partial f(\bar{x})$ when taking the limit $x^k\to\bar{x}$ in \eqref{eq:1.7}. Further, the feasibility of $\bar{x}$ follows from the continuity of $g$ such that we immediately obtain that $\bar{x}$ is a KKT point with multiplier $u\coloneqq 0$.\\
For the case $I\neq \emptyset$, we first show that the sequence $\{ u^k \}_{k=1}^\infty$ is bounded. 
To this end, let us assume the contrary. Then, without loss of generality, $\lVert u^k \rVert \to \infty$ as $k\to\infty$. Now, consider the bounded sequence $\{ \bar{u}^k \}_{k=1}^\infty$ given by
\vspace{-0.2em}
\[  
\bar{u}^k \coloneqq \frac{u^{k}}{\lVert u^k \rVert}. 
\]
Due to the boundedness, we may assume without loss of generality that $\bar{u}^k \to \bar{u}$ with $\lVert \bar{u} \rVert = 1$ and $\bar{u}\geq 0$.
Since $f$ is locally Lipschitz continuous and ${x^{k}}\to\bar{x}$, we have by parts \ref{propa} and \ref{propc} of Proposition~\ref{prop:eig_clarke_subdiff} that any sequence $\{ s^k \}_{k=1}^\infty$ with $s^k \in \partial f(x^k)$ is bounded, which implies
\vspace{-0.5em}
\[ 
\frac{s^k}{\lVert u^k \rVert} \to 0 \quad\text{as}\,\,\, k\to\infty. 
\]
Dividing by $\lVert u^k \rVert$ and passing to the limit in \eqref{eq:1.7} therefore yields
\[
0 \in\, \sum_{j\in {I}} \bar{u}_j \partial g_j(\bar{x}), 
\vspace{-0.2em} 
\]
where we used  $\epsilon^k \to 0$ and Proposition \ref{prop:eig_clarke_subdiff}\ref{propb} as $g$ is locally Lipschitz continuous. 
This contradicts subMFC, and hence, the sequence $\{ u^k \}_{k=1}^\infty$ has to be bounded. Thus, without loss of generality, $u^k \to \widetilde{u}$. In turn, as $k\to\infty$, it follows from \eqref{eq:akkt2}, \eqref{eq:1.7}, Proposition \ref{prop:eig_clarke_subdiff}, $\delta^k \to 0$ and the fact that ${I} \subseteq I(\bar{x})$ by the continuity of $g$, that $\bar{x}$ is a KKT point of \eqref{P} with multiplier~$\widetilde{u}$.
\end{proof}

The subsequent example shows that the implication in Theorem \ref{thm:bqc_implies_kkt} cannot be inverted, proving that subMFC is not simply a reformulation of the KKT conditions. Moreover, the example reveals that subMFC is not a necessary optimality condition. Rather, although part \ref{subMFCII} of subMFC is a sequential necessary optimality condition that also holds at any KKT point, subMFC as a whole is not. For this reason, it also does not belong to the class of sequential optimality conditions in the classical sense \cite{RA3,RA4,RA7,EB1}. Thus, the classification of subMFC as problem qualification is justified.

\begin{example}
\label{ex:sub_mfc_kkt}
We consider the Lipschitz continuous problem
\begin{align*}
\min\, 3x_1-2x_2\quad 
\text{ s.t. }\quad& g_1(x_1,x_2) \coloneqq |x_1|-x_2 \leq 0, \\
& g_2(x_1,x_2) \coloneqq -|x_1|+x_2 \leq 0, \\
& g_3(x_1,x_2) \coloneqq x_1^2+(x_2+1)^2-1 \leq 0, \\
& g_4(x_1,x_2) \coloneqq -x_1^2-(x_2+1)^2+1 \leq 0
\end{align*}
with its only feasible point $\bar{x}=(\bar{x}_1,\bar{x}_2)=(0,0)$. First, we see that $\bar{x}$ is a KKT point (and thus also an AKKT point), as
\[
\begin{pmatrix}0 \\ 0\end{pmatrix} \\
\in \begin{pmatrix} 3 \\ -2 \end{pmatrix} + u_1  \begin{pmatrix} [-1,1] \\ -1 \end{pmatrix} + u_2 \begin{pmatrix} [-1,1] \\ 1 \end{pmatrix} + u_3 \begin{pmatrix} 0 \\ 2 \end{pmatrix} + u_4 \begin{pmatrix} 0 \\ -2 \end{pmatrix} 
\]
can be fulfilled by choosing $u=(0,3,0,\frac{1}{2})$ and the subgradient $(-1,1)\in \partial g_2(\bar{x})$.
To show that subMFC does not hold at $\bar{x}$, let us pick a sequence $\{(x^k,u^k,\delta^k,\epsilon^k)\}_{k=1}^\infty$ satisfying \ref{subMFCII} of subMFC. Then, $(x^k, \delta^k,\epsilon^k) \to (\bar{x},0,0)$ holds for $k\to\infty$ and for all $k\in\mathbb{N}$ we have
\begin{gather*}
\hspace{-0.802em}
\begin{pmatrix}\epsilon^k_1 \\ \epsilon^k_2\end{pmatrix}
\in \begin{pmatrix} 3 \\ -2 \end{pmatrix} 
+ u^k_1  \begin{pmatrix} \partial(|x_1^k|) \\ -1 \end{pmatrix} 
+ u^k_2 \begin{pmatrix} -\partial(|x_1^k|) \\ 1 \end{pmatrix} 
+ u^k_3 \begin{pmatrix} 2x_1^k \\ 2(x_2^k+1) \end{pmatrix} 
+ u^k_4 \begin{pmatrix} -2x_1^k \\ -2(x_2^k+1) \end{pmatrix}, 
\numberthis \label{eq:ex_acq} 
\\
u^k_i \geq 0, \quad g_i(x_1^k,x_2^k) \leq \delta^k_i,\quad u^k_i (g_i(x_1^k,x_2^k) - \delta^k_i) = 0, \quad \forall i=1,\dotsc,4,
\vspace{-0.3em}
\end{gather*}
where we use $\partial(|x_1^k|)$ as an abbreviation for $\partial|\cdot|(x_1^k)$, i.e., for the subdifferential of the function $x\mapsto |x|$ at $x_1^k$.
We now want to find out for which sets $I\subseteq I(\bar{x})=\{1,2,3,4\}$ part~\ref{subMFCI} of subMFC can be satisfied.
As we have 
\[ 
\partial g_1(\bar{x}) = \begin{pmatrix} [-1,1] \\ -1 \end{pmatrix}, \partial g_2(\bar{x}) = \begin{pmatrix} [-1,1] \\ 1 \end{pmatrix}, \nabla g_3(\bar{x}) = \begin{pmatrix} 0 \\ 2 \end{pmatrix}, \nabla g_4(\bar{x}) = \begin{pmatrix} 0 \\ -2 \end{pmatrix},  
\]
part \ref{subMFCI} of subMFC clearly only holds if 
\[ 
\text{(a)}\,\,  I(x^k,\delta^k) \subseteq \{ 1,4 \} \quad\text{or}\quad \text{(b)}\,\, I(x^k,\delta^k) \subseteq \{ 2,3 \}
\]
is valid. This implies for $u^k=(u_1^k,u_2^k,u_3^k,u_4^k)$ that 
\[ 
\text{(a)}\,\, u_1^k\geq 0, u_2^k=0, u_3^k=0, u_4^k\geq 0 \quad\text{or}\quad \text{(b)}\,\, u_1^k=0,u_2^k\geq 0,u_3^k\geq 0, u_4^k = 0. 
\]
In case (a), we immediately see that the second line of \eqref{eq:ex_acq} reads as
$
\epsilon^k_2 = -2 - u^k_1 -2(x_2^k+1) u^k_4, 
$
and hence, $\epsilon^k_2\to 0$ as $k\to\infty$ cannot be satisfied. Therefore, only case (b) has to be considered in the following, meaning that \eqref{eq:ex_acq} reads as 
\[
\begin{pmatrix}\epsilon^k_1 \\ \epsilon^k_2\end{pmatrix} \\
\in \begin{pmatrix} 3 \\ -2 \end{pmatrix} + u^k_2 \begin{pmatrix} -\partial(|x_1^k|) \\ 1 \end{pmatrix} + u^k_3 \begin{pmatrix} 2x_1^k \\ 2(x_2^k+1) \end{pmatrix}. \numberthis \label{eq:ex_acq1}
\]
As $x^k \to \bar{x}=(0,0)$ for $k\to\infty$, we may assume without loss of generality (by passing to a subsequence if necessary) that one of the three cases $x_1^k >0$, $x_1^k=0$ or $x_1^k<0$ holds for all $k\in\mathbb{N}$.
\begin{itemize}
\item If $x_1^k > 0$, system \eqref{eq:ex_acq1} reads as
\[
\begin{pmatrix}\epsilon^k_1 \\ \epsilon^k_2\end{pmatrix} \\
= \begin{pmatrix} 3 \\ -2 \end{pmatrix} + u^k_2 \begin{pmatrix} -1 \\ 1 \end{pmatrix} + u^k_3 \begin{pmatrix} 2x_1^k \\ 2(x_2^k+1) \end{pmatrix},
\]
leading to
\[ 
\epsilon^k_2 = -2 + u^k_2 + 2(x_2^k+1)u^k_3 = -2 + 3 - \epsilon^k_1 + 2x_1^k u^k_3 + 2(x_2^k+1)u^k_3. 
\]
The right-hand side is clearly bounded away from 0 for large $k\in\mathbb{N}$ and hence, $\epsilon^k_2\to 0$ as $k\to\infty$ cannot be satisfied.
\item If $x_1^k = 0$, equation \eqref{eq:ex_acq1} amounts to
\[
\begin{pmatrix}\epsilon_1^k \\ \epsilon^k_2\end{pmatrix} \\
\in \begin{pmatrix} 3 \\ -2 \end{pmatrix} + u^k_2 \begin{pmatrix} [-1,1] \\ 1 \end{pmatrix} + u^k_3 \begin{pmatrix} 0 \\ 2(x_2^k+1) \end{pmatrix}, \numberthis\label{eq:ex_acq2}
\]
which yields
$\epsilon^k_1 = 3 + u_2^k \cdot s$ for some $s\in [-1,1]$.
To guarantee $\epsilon^k_1\to 0$ as $k\to\infty$, only $s\in [-1,0)$ is possible. Thus, we obtain
$u_2^k = (\epsilon^k_1-3)s^{-1}$ for some $s\in [-1,0)$.
Plugging the latter into the second line of \eqref{eq:ex_acq2} gives
\[ 
\epsilon^k_2 = -2 + u^k_2 + 2(x_2^k+1)u^k_3 = -2 + (\epsilon^k_1-3)s^{-1} + 2(x_2^k+1)u^k_3 \geq -2 + 3-\epsilon^k_1 + 2(x_2^k+1)u^k_3.
\]
This again is bounded away from $0$ for large $k\in\mathbb{N}$, such that $\epsilon^k_2\to 0$ as $k\to\infty$ cannot hold.
\item In the case $x_1^k < 0$, equation \eqref{eq:ex_acq1} reads as
\[
\begin{pmatrix}\epsilon^k_1 \\ \epsilon^k_2\end{pmatrix} \\
= \begin{pmatrix} 3 \\ -2 \end{pmatrix} + u^k_2 \begin{pmatrix} 1 \\ 1 \end{pmatrix} + u^k_3 \begin{pmatrix} 2x_1^k \\ 2(x_2^k+1) \end{pmatrix},
\]
which yields $\epsilon^k_2 = -2 + u^k_2 + 2(x_2^k+1)u^k_3$ and thus
\begin{align*}
\hspace{-0.1em}
u^k_3 = \frac{2 + \epsilon^k_2 - u^k_2}{2(x_2^k+1)},\,
 \epsilon^k_1 
&= 3 + u^k_2 + x_1^k \,\frac{2 + \epsilon^k_2 - u^k_2}{x_2^k+1} 
= 3+\left(1-\frac{x_1^k}{x_2^k+1}\right)u_2^k + \frac{x_1^k(2+\epsilon_2^k)}{x_2^k+1}.
\end{align*} 
Clearly, it is again not possible to satisfy $\epsilon^k_1\to 0$ as $k\to\infty$.
\end{itemize}
Thus, it is impossible to satisfy both \ref{subMFCI} and \ref{subMFCII} of subMFC at the same time. 
Altogether, this means that subMFC cannot hold at the KKT point $\bar{x}$.
\end{example}

\begin{remark}
As Theorem \ref{thm:bqc_implies_kkt} holds for any $\bar{x}$ that is feasible for \eqref{P}, we can state analogous results for $\bar{x}$ being a local minimizer or for $\bar{x}$ being an AKKT point of \eqref{P}. For these special cases, the theorem reveals that subMFC can be used as a condition to guarantee that a local minimizer (resp. an AKKT point) fulfills the KKT conditions. This essentially means that subMFC still shares the defining property of a CQ (resp. of a strict CQ), even though it is not a CQ itself. 
\end{remark}

\begin{remark}
Of course, we could also state Definitions \ref{def:akkt} and \ref{def:subMFC} with respect to other subdifferentials than Clarke's subdifferential, such as the regular subdifferential $\widehat{\partial}(\cdot)$ or the limiting subdifferential $\partial^L(\cdot)$. As these subdifferentials are generally smaller than Clarke's subdifferential, the definition of an AKKT point would become more restrictive. For subMFC, the positive linear independence required for part \ref{subMFCI} would generally be easier to fulfill, whereas part \ref{subMFCII} would become harder to fulfill. Hence, it is not straightforward to derive a general rule whether subMFC becomes weaker or stronger when using a different subdifferential.
\end{remark}

To end this section, we want to illustrate how subMFC can be applied. Simultaneously, this part also clarifies how subMFC can be used to treat problems with equality constraints. For this purpose, we first take a look at a standard equality constrained example in which MFCQ fails to demonstrate the benefit of using subMFC.

\newpage

\begin{example}
\label{ex:equality1}
Consider an optimization problem with the feasible set given by the equality constraint $x^2=0$, which can equivalently be written as feasible set of an instance of optimization problem~\eqref{P} by splitting the equation into two inequalities:
\begin{align*}
\min \,  x^2 \quad \text{ s.t. } \quad &g_1(x) \coloneqq x^2 \leq 0,\,\, g_2(x) \coloneqq -x^2 \leq 0.
\end{align*}
In the only feasible point and global minimizer $\bar{x}=0$, we obtain $I(\bar{x})=\{1,2\}$. Clearly, MFCQ does not hold at $\bar{x}$. Concerning subMFC, we first take a look at the AKKT conditions in $\bar{x}$ for the above problem, which require $(x^k, \delta^k,\epsilon^k) \to (\bar{x},0,0)$ for $k\to\infty$ as well as
\begin{gather*}
\epsilon^k = 2x^k + u_1^k (2x^k) + u_2^k (-2x^k), \\
\begin{alignat*}{4}
u^k_1 \geq 0,&& (x^k)^2 \leq \delta^k_1,&& u^k_1 ((x^k)^2 - \delta^k_1) = 0, \phantom{aaaaaaaaaaaa}\\
u^k_2 \geq 0,&&\quad -(x^k)^2 \leq \delta^k_2,&&\quad u^k_2 (-(x^k)^2 - \delta^k_2) = 0 \phantom{,aaaaaaaaaaaa}
\end{alignat*}
\end{gather*}
for all $k\in\mathbb{N}$.
To satisfy subMFC at $\bar{x}$, we need to find a sequence $\{(x^k,u^k,\delta^k,\epsilon^k)\}_{k=1}^\infty$ that fulfills the above conditions such that either $I^g(x^k,\delta^k) = \emptyset$ or 
\[ 
0 \displaystyle \neq \sum_{j\in I^g(x^k,\delta^k)} v_j \nabla g_j(\bar{x})
\]
holds for all $v\geq 0$, $v\neq 0$.
Due to $\nabla g_1(\bar{x})=0$ and $\nabla g_2(\bar{x}) = 0$, this implies that we need $I^g(x^k,\delta^k) = \emptyset$ for subMFC to hold at $\bar{x}$. Indeed, a sequence fulfilling both the AKKT conditions and $I^g(x^k,\delta^k) = \emptyset$ can be found, e.g.
\[ 
\{(x^k,u^k,\delta^k,\epsilon^k)\}_{k=1}^\infty \coloneqq \left\{0,\left(0,0\right),\left(\frac{1}{k},\frac{1}{k}\right),0\right\}_{k=1}^\infty. 
\]
\end{example}

Thus, the example shows that we can exclude unwanted gradients from the index set $I^g(x^k,\delta^k)$ by setting $\delta^k$ accordingly. The same idea can also be applied if the point of interest is not the unrestricted minimizer of the objective function, i.e., when $I^g(x^k,\delta^k) \neq \emptyset$ is required, as the following example shows.

\begin{example}
\label{ex:two_eqs}
Consider an optimization problem with the objective function $f(x_1,x_2) \coloneqq x_1-2x_2$ and the feasible set given by the equality constraints $h_1(x)=h_1(x_1,x_2)\coloneqq |x_1|-x_2=0$ and $h_2(x)=h_2(x_1,x_2)\coloneqq x_1^2+(x_2+1)^2-1=0$, for which the only feasible point is $\bar{x}=(\bar{x}_1,\bar{x}_2)=(0,0)$. We obtain the (sub)gradients
\[ 
\partial h_1(x)=  \begin{pmatrix} \partial(|x_1|) \\ -1 \end{pmatrix}
\,\,\text{with}\,\,  \partial(|x_1|)  = \begin{cases}
    -1, & \text{for } x_1<0, \\
    [-1,1], & \text{for } x_1 = 0, \\
    1, & \text{for } x_1>0
  \end{cases}
\,\,\text{ and }\,\, \nabla h_2(x)= \begin{pmatrix} 2x_1 \\ 2(x_2+1) \end{pmatrix},
\]
amounting to
\[ 
\partial h_1(\bar{x})= \begin{pmatrix}  [-1,1] \\ -1 \end{pmatrix} \,\,\text{ and }\,\,  \nabla h_2(\bar{x})= \begin{pmatrix} 0 \\ 2 \end{pmatrix}
\]
at $\bar{x}$. Since the (sub)gradients $(0,-1) \in\partial h_1(\bar x)$ and $\nabla h_2(\bar{x})$ 
are linearly dependent, LICQ does not hold at $\bar{x}$. In turn, even MFCQ tailored to problems including equality constraints is violated at $\bar{x}$. \\ The reformulation of the problem into an instance of \eqref{P}, obtained by splitting the equations into two inequalities, is given by the constraint set in Example \ref{ex:sub_mfc_kkt} together with the objective function $f(x_1,x_2) = x_1-2x_2$.
To verify subMFC at $\bar{x}$, we write down the AKKT conditions at $\bar{x}$, which consist of $(x^k, \delta^k,\epsilon^k) \to (\bar{x},0,0)$ as $k\to\infty$ and 
\begin{gather*}
\begin{pmatrix}\epsilon^k_1 \\ \epsilon^k_2\end{pmatrix}
\in \begin{pmatrix} 1 \\ -2 \end{pmatrix} + u^k_1  \begin{pmatrix} \partial(|x_1^k|) \\ -1 \end{pmatrix} + u^k_2 \begin{pmatrix} -\partial(|x_1^k|) \\ 1 \end{pmatrix} + u^k_3 \begin{pmatrix} 2x_1^k \\ 2(x_2^k+1) \end{pmatrix} + u^k_4 \begin{pmatrix} -2x_1^k \\ -2(x_2^k+1) \end{pmatrix}, \\
u^k_i \geq 0, \quad g_i(x_1^k,x_2^k) \leq \delta^k_i,\quad u^k_i (g_i(x_1^k,x_2^k) - \delta^k_i) = 0, \quad \forall i=1,\dotsc,4
\end{gather*}
for all $k\in\mathbb{N}$.
It can be checked that one sequence to fulfill this system is given by
\[ 
\{({x}^k,u^k, \delta^k, \epsilon^k)\}_{k=1}^\infty \coloneqq \left\{\left(\left(\frac{1}{k},0\right),\left(0,1,\frac{1}{2},0\right),\left(\frac{2}{k},-\frac{1}{k},\frac{1}{k^2},0\right),\left( \frac{1}{k} ,0\right)\right)\right\}_{k=1}^\infty. 
\]
For this sequence we have $I(x^k,\delta^k) = \{ 2,3 \}$ for all $k\in\mathbb{N}$. Thus, setting $I\coloneqq \{ 2,3 \}$, we see that
\[ 
0 \displaystyle \notin \sum_{j\in {I}} v_j \partial g_j(\bar{x}) = v_2 \begin{pmatrix} -\partial(|\bar{x}_1|) \\ 1 \end{pmatrix} + v_3 \begin{pmatrix} 2\bar{x}_1 \\ 2(\bar{x}_2+1) \end{pmatrix} = v_2 \begin{pmatrix} [-1,1] \\ 1 \end{pmatrix} + v_3 \begin{pmatrix} 0 \\ 2 \end{pmatrix} 
\]
is satisfied for all $(v_2,v_3) \neq (0,0)$ with $v_2\geq 0,v_3 \geq 0$, i.e., subMFC holds at $\bar{x}$.
\end{example}

In summary, the constraints that are active in $\bar{x}$ have linearly dependent (sub)gradients, which leads to a violation of MFCQ, no matter whether the constraints are posed as equality constraints or as split inequality constraints. On the opposite, the split into inequality constraints allows subMFC to choose an appropriate selection of the (sub)gradients. Precisely, we only pick those linearly independent (sub)gradients which are crucial for fulfilling the AKKT conditions, while all others are neglected by choosing the corresponding component of $\delta^k$ accordingly. In particular, we have seen that subMFC is able to identify KKT points starting from AKKT points in situations where standard CQs are not readily available. The following section is dedicated to showing that this benefit also occurs in situations where far weaker CQs and problem qualifications are violated.

\section{Comparison to constraint qualifications and calmness}
\label{sec:other_conds}
In order to show that subMFC is a reasonable and powerful condition, the next aim is to prove in Sections \ref{subsec:mfcq}--\ref{subsec:quasinorm} that it is automatically fulfilled for a nontrivial class of optimization problems, as many well-known CQs strictly imply subMFC at local minimizers. Thereafter, we demonstrate in Sections \ref{subsec:errorbound}--\ref{subsec:calmness} that several weak conditions are not implied by subMFC or known to be independent of subMFC. For a clearer presentation of this part, the corresponding detailed calculations are provided in the appendix. Finally, in Section \ref{subsec:summarygraph}, a summary of the investigated dependencies is given.

\subsection{Mangasarian-Fromovitz constraint qualification}
\label{subsec:mfcq}
Recall the definition of MFCQ, Definition \ref{def:mfcq}. As indicated before, subMFC can be interpreted as a weakened MFCQ which aligns with the next result.

\begin{theorem}
\label{thm:mfcq_subMFC}
Let $\bar{x}$ be a local minimizer of \eqref{P}. Then, MFCQ at $\bar{x}$ strictly implies subMFC at~$\bar{x}$.
\end{theorem}

\begin{proof}
As $\bar{x}$ is an AKKT point by Theorem \ref{thm:locmin_akkt}, there exists a sequence $\{(x^k,u^k,\delta^k,\epsilon^k)\}_{k=1}^\infty$ that fulfills system \eqref{eq:akkt1}--\eqref{eq:akkt5}. Recalling Remark \ref{rem:seq_const}, we can take $I \coloneqq I(x^k,\delta^k)$ to fulfill part~\ref{subMFCII} of subMFC at $\bar{x}$. 
By the continuity of $g$ and due to $\delta^k \to 0$, it follows that $I\subseteq I(\bar{x})$.
Thus, if $I(\bar{x})=\emptyset$, then $I=\emptyset$ such that part \ref{subMFCI} of subMFC is fulfilled automatically.
Otherwise, if $I(\bar{x})\neq\emptyset$, then by MFCQ we have for all $v\geq 0, v \neq 0$ that
\[ 
0 \displaystyle \notin \sum_{j\in I(\bar{x})} v_j \partial g_j(\bar{x}) \overset{I \subseteq I(\bar{x})}{\implies} 0 \displaystyle \notin \sum_{j\in I} v_j \partial g_j(\bar{x}). 
\]
Hence, part \ref{subMFCI} of subMFC holds at $\bar{x}$ as well. Finally, it follows from Example \ref{ex:equality1} that this implication is strict.
\end{proof}

Related to MFCQ, a well-known problem qualification in standard nonlinear programming is the strict Mangasarian-Fromovitz CQ (SMFCQ) which, given a KKT point, ensures uniqueness of the associated Lagrange multipliers \cite{JK1}. 
Clearly, as SMFCQ presumes the existence of a Lagrange multiplier, it depends on the objective function and, thus, is not a CQ but a problem qualification. To emphasize this circumstance, it was also named strict Mangasarian-Fromovitz condition (SMFC) later on \cite{GW1}.
Since this problem qualification is weaker than LICQ but stronger than MFCQ \cite{JK1}, it trivially implies subMFC as well.

\subsection{Constant positive linear dependence constraint qualification}
\label{subsec:cpld}
We can adapt the more general CPLD definition from \cite[Definition 1.1]{MX1} for nonsmooth problems to our setting in \eqref{P} as follows.

\begin{definition}
\label{def:CPLD}
\textit{CPLD} is said to hold at a point $\bar{x}$ feasible for \eqref{P} if, for any nonempty subset $J \subseteq I(\bar{x})$ and any choice of subgradients $s_j\in \partial g_j(\bar{x}), j\in J,$ such that there exists $v \geq 0, v\neq 0$ with
\[ 
0 = \sum_{j\in J} v_j s_j, 
\]
then, for all $k\in\mathbb{N}$ sufficiently large and for all sequences $\{ x^k \}_{k=1}^\infty$, $\{ s^k \}_{k=1}^\infty$ satisfying $x^k \neq \bar{x}$, $x^k \to \bar{x}$ as $k\to\infty$, $s^k_j \in \partial g_j(x^k)$ and $s^k_j \to s_j$ as $k\to\infty$, there exists ${v^k}\neq 0$ such that 
\[ 
0 = \sum_{j\in J} v_j^k s_j^k. 
\]
\end{definition}

It is known that CPLD is a CQ \cite{MX1} and, by definition, it is straightforward to see that MFCQ implies CPLD. In \cite[Proposition 4.1]{MX1}, it is further shown that CPLD is also implied by CRCQ. Both of these implications are strict \cite{RA1}. To prove that CPLD implies subMFC at local minimizers, we make use of the following auxiliary result \cite[Lemma 3.1]{EB1} that can be viewed as a corollary of Carathéodory's Lemma.

\begin{lemma}
\label{lem:lin_indep}
Assume that
$v = \sum_{j=1}^{q} \lambda_j v_j$
with $v_j\in\mathbb{R}^n$ and $\lambda_j \geq 0$ for all $j=1,\dotsc,q$. Then, there exist $J\subseteq\{1,\dots,q\}$ and $\bar{\lambda}_j\geq 0$ for all $j \in J$ such that:
\begin{enumerate}[label=(\alph*)]
\item $v = \sum_{j \in J} \bar{\lambda}_j v_j$;
\item The vectors $\{ v_j \}_{j \in J}$ are linearly independent.
\end{enumerate}
\end{lemma}

The next theorem shows that CPLD implies subMFC at a local minimizer if the subdifferentials of the active constraints are singletons in that local minimizer.

\begin{theorem}
\label{thm:cpld_subMFC}
Let $\bar{x}$ be a local minimizer of \eqref{P}.
If $\partial g_j(\bar{x})$ is a singleton for all $j\in I(\bar{x})$, then CPLD at $\bar{x}$ strictly implies subMFC at $\bar{x}$.
\end{theorem}

\begin{proof}
Since $\bar{x}$ is an AKKT point by Theorem \ref{thm:locmin_akkt}, there is a sequence $\{(x^k,u^k, \delta^k, \epsilon^k)\}_{k=1}^\infty$ that fulfills \eqref{eq:akkt1}--\eqref{eq:akkt5}.
This yields in particular
\begin{align*}
\epsilon^k = t^k + \sum_{j\in I(x^k,\delta^k)}  u^k_j s^k_j \qquad\text{for some }\, t^k \in \partial f(x^k),\,\, s^k_j \in \partial g_j(x^k).
\end{align*} 
By Lemma \ref{lem:lin_indep} and Remark \ref{rem:seq_const}, we may assume that there exist an index set $J$, an infinite index set $N\subseteq\mathbb{N}$ such that $J\subseteq I(x^k,\delta^k)$ for all $k\in N$, and multipliers $\bar{u}^k\geq 0$ for all $k\in N$ such that:
\begin{enumerate}[label=(\alph*)]
\item $\epsilon^k = t^k + \sum_{j\in J} \bar{u}^k_j s^k_j$;
\item The vectors $\{ s^k_j \}_{j\in J}$ are linearly independent.\vspace{0.5em}
\end{enumerate}
Suppose now that subMFC does not hold at $\bar{x}$. Then, there exists $v \geq 0, v \neq 0$ such that
\begin{align*}
0 = \sum_{j\in I(x^k,\delta^k)} v_j s_j,
\end{align*} 
where $\partial g_j(\bar{x}) = \{ s_j \}$ for all $j\in I(x^k,\delta^k)$ by assumption, as we have $I(x^k,\delta^k)\subseteq I(\bar{x})$ by $\delta^k \to 0$ and the continuity of $g$. Moreover, by parts \ref{propa} and \ref{propc} of Proposition \ref{prop:eig_clarke_subdiff}, we find that, for all $j\in I(x^k,\delta^k)$, the sequence $\{ s^k_j \}_{k=1}^\infty$ is convergent for $k\to\infty$ without loss of generality. Due to Proposition \ref{prop:eig_clarke_subdiff}\ref{propb}, $x^k\to\bar{x}$ and $\partial g_j(\bar{x}) = \{ s_j \}$, this particularly means that $s^k_j \to s_j$ for $k\to\infty$. As CPLD holds at $\bar{x}$, the vectors $\{ s_j^k \}_{j\in I(x^k,\delta^k)}$ must therefore be linearly dependent.\\
If $J = I(x^k,\delta^k)$, this is a contradiction to (b), and hence, subMFC holds at $\bar{x}$.
Otherwise, we have $I(x^k,\delta^k) \setminus J \neq \emptyset$.
As this set contains the indices of $I(x^k,\delta^k)$ that are not used in representation (a), we set
\begin{align*}
(\widetilde{u}^k_j, \widetilde{\delta}^k_j) \coloneqq 
\begin{cases}
	( 0, c\delta^k_j) & \forall j \in I(x^k,\delta^k) \setminus J \text{ with } \delta^k_j > 0,\\
 	( 0, -c\delta^k_j) & \forall j \in I(x^k,\delta^k) \setminus J \text{ with } \delta^k_j < 0,\\
 	( 0, \frac{1}{k}) & \forall j \in I(x^k,\delta^k) \setminus J \text{ with } \delta^k_j = 0,\\
 	( \bar{u}^k_j, \delta^k_j) & \forall j \in J,\\
 	( 0, \delta^k_j) &\forall j \notin I(x^k,\delta^k)
\end{cases}
\end{align*}
for some $c>1$.
Then, the sequence $\{(x^k,\widetilde{u}^k, \widetilde{\delta}^k,\epsilon^k)\}_{k=1}^\infty $ fulfills
\vspace{-0.6em}
\begin{gather}
 \epsilon^k \in \partial f(x^k) + \sum_{j=1}^q  \widetilde{u}^k_j \partial g_j(x^k), \label{eq:1.8a}\\
 g(x^k) \leq \widetilde{\delta}^k, \quad \widetilde{u}^k \geq 0,\quad (\widetilde{u}^k)^\top (g(x^k)-\widetilde{\delta}^k) = 0, \vspace{0.5em} \label{eq:1.8}\\
 (x^k, \widetilde{\delta}^k,\epsilon^k) \to (\bar{x},0,0) \quad \text{as}\,\,\, k\to\infty. \label{eq:1.11}
\end{gather}
Here, \eqref{eq:1.8a} follows from (a). Additionally, \eqref{eq:1.8} follows from $\delta^k_j \leq \widetilde{\delta}^k_j$ and the fact that $\widetilde{u}^k_j > 0$ is possible only for $j\in J\subseteq I(x^k,\delta^k)$ for which we have $g_j(x^k)-\widetilde{\delta}^k_j = g_j(x^k)- \delta^k_j=0$.
For the new perturbation vector $\widetilde{\delta}^k$ defined above, it holds that
\vspace{-0.2em}
\begin{alignat*}{3}
g_j(x^k) - \widetilde{\delta}^k_j &< g_j(x^k)- \delta^k_j = 0 \qquad && \forall j \in I(x^k,\delta^k) \setminus J, \\
g_j(x^k) - \widetilde{\delta}^k_j &= g_j(x^k)- \delta^k_j = 0 && \forall j \in J,  \\
g_j(x^k) - \widetilde{\delta}^k_j &= g_j(x^k)- \delta^k_j < 0 && \forall j \notin I(x^k,\delta^k), 
\end{alignat*}
and hence, $I(x^k,\widetilde{\delta}^k) = J$. The beginning of the proof can now be repeated for the new sequence $\{(x^k,\widetilde{u}^k,$ $\widetilde{\delta}^k,\epsilon^k)\}_{k=1}^\infty $ as \eqref{eq:1.8a}--\eqref{eq:1.11} corresponds to \eqref{eq:akkt1}--\eqref{eq:akkt5} being fulfilled with $\widetilde{u}^k$ and $\widetilde{\delta}^k$ instead of $u^k$ and $\delta^k$. As a consequence, we obtain that subMFC holds at $\bar{x}$.\\
Finally, we demonstrate that this implication is strict. To this end, consider the optimization problem from Example \ref{ex:equality1}, for which it is already known that subMFC holds at the global minimizer $\bar{x}=0$. Thus, we now verify that CPLD is violated at $\bar{x}$. According to Definition~\ref{def:CPLD}, we have to look at any nonempty subset $J \subseteq I(\bar{x}) = \{ 1,2 \}$. For the choice $J \coloneqq \{ 1 \}$, we find $ \nabla g_1(\bar{x})=0 $, and hence, there exists $v_1>0$ with $0 = v_1 \cdot \nabla g_1(\bar{x})$. The sequences $\{ x^k \}_{k=1}^\infty$ and $\{ s^k \}_{k=1}^\infty$ given by $x^k \coloneqq \frac{1}{k}$ and $s^k \coloneqq \frac{2}{k}$ satisfy $x^k \neq 0$, $x^k \to 0$ as $k\to\infty$, $s^k = \nabla g_1(x^k)$ and $s^k \to s$ as $k\to\infty$. However, there exists no $v_1^k \neq 0$ such that for all $k\in\mathbb{N}$ sufficiently large we have 
$ 0 = v_1^k \cdot s^k $. Hence, CPLD does not hold. 
\end{proof}

Note that the assumption that $\partial g_j(\bar{x})$ is a singleton for all $j\in I(\bar{x})$ is required as otherwise it would not be guaranteed that the subgradients required for the AKKT conditions and the subgradients required in the assumed contradiction to subMFC coincide. However, taking a closer look at the proof, we can relax this assumption slightly.

\begin{corollary}
\label{cor:cpld_subMFC}
Let $\bar{x}$ be a local minimizer of \eqref{P}. Let $\{(x^k,u^k, \delta^k, \epsilon^k)\}_{k=1}^\infty$ be a sequence such that Definition \ref{def:akkt} is fulfilled at $\bar{x}$ with $I(x^k,\delta^k)$ being constant for all $k\in\mathbb{N}$.
If $\partial g_j(\bar{x})$ is a singleton for all $j\in I(x^k,\delta^k)$, then CPLD at $\bar{x}$ strictly implies subMFC at $\bar{x}$.
\end{corollary}

As $I(x^k,\delta^k) \subseteq I(\bar{x})$ holds due to the continuity of $g$ and $\delta^k \to 0$, Corollary \ref{cor:cpld_subMFC} generally requires less subdifferentials to be a singleton than Theorem \ref{thm:cpld_subMFC}. Moreover, the additional assumptions concerning the concrete sequence are not restrictive since they can always be fulfilled, see Theorem~\ref{thm:locmin_akkt} and Remark~\ref{rem:seq_const}.

\subsection{Pseudonormality and quasinormality}
\label{subsec:quasinorm}
Next, we consider the well-known definitions of pseu\-donormality and quasinormality, introduced in \cite{DB1} and \cite{MH1}. They can be extended to the Lipschitz continuous setting in \eqref{P} as follows.

\begin{definition}
\label{def:pq_normality}
A point $\bar{x}$ feasible for \eqref{P} is said to satisfy
\begin{itemize}
\item[(a)] \textit{pseudonormality} if there exists no multiplier $v \geq 0, v\neq 0$ such that:
\vspace{0.5em}
\begin{itemize}
\item[(a1)] $ 0 \displaystyle \in \sum_{j\in I(\bar{x})} v_j \partial g_j(\bar{x})$;
\item[(a2)] There exists a sequence $\{ x^k \}_{k=1}^\infty$ with $x^k \to \bar{x}$ as $k\to\infty$ such that
\[ 
\sum_{j\in I(\bar{x})} v_j g_j({x}^k) >0 
\]
holds for all $k\in\mathbb{N}$.
\vspace{0.5em}
\end{itemize}
\newpage
\item[(b)] \textit{quasinormality} if there exists no multiplier $v \geq 0, v\neq 0$ such that:
\vspace{0.5em}
\begin{enumerate}[label=(b\arabic*)]
\item \label{item:quasinorm1} $ 0 \displaystyle \in \sum_{j\in I(\bar{x})} v_j \partial g_j(\bar{x})$;
\item \label{item:quasinorm2} There exists a sequence $\{ x^k \}_{k=1}^\infty$ with $x^k \to \bar{x}$ as $k\to\infty$ such that
\begin{align*}
v_j g_j({x}^k) &> 0 \quad\forall j \text{ with } v_j > 0
\end{align*}
holds for all $k\in\mathbb{N}$.
\end{enumerate}
\end{itemize}
\end{definition}

The definitions immediately show that both conditions are weaker than MFCQ. Moreover, they show that pseudonormality implies quasinormality, where the reverse statement is not true by \cite[Example 3.1]{DB1}. If all constraint functions $g_j, j=1,\dots,q$, of \eqref{P} are concave, then any feasible point of \eqref{P} satisfies pseudonormality \cite[Proposition 3]{JY11}. Further, in~\cite{RA1} it is shown that CPLD neither implies nor is implied by pseudonormality. In contrast, quasinormality is implied by CPLD for the smooth case \cite[Theorem 3.1]{RA1}. However, this implication does not hold anymore for the nonsmooth case \cite{R1}, which is the reason why we considered CPLD separately in Section \ref{subsec:cpld}. Fortunately, it turns out that even quasinormality implies subMFC.

\begin{theorem}
\label{thm:quasinorm}
Let $\bar{x}$ be a local minimizer of \eqref{P}. Then, quasinormality at $\bar{x}$ strictly implies subMFC at $\bar{x}$.
\end{theorem}

\begin{proof}
From Remarks \ref{rem:delta_greater_0} and \ref{rem:seq_const} it follows that there exists a sequence $\{(x^k,u^k,\delta^k,\epsilon^k)\}_{k=1}^\infty$ which fulfills system \eqref{eq:akkt1}--\eqref{eq:akkt5} with $\delta^k_j>0$ for all $j=1,\dotsc,q$ and $I(x^k,\delta^k)$ being constant for all $k\in\mathbb{N}$. Suppose that $\bar{x}$ satisfies quasinormality, while subMFC does not hold at $\bar{x}$. For $I \coloneqq I(x^k,\delta^k)$, it follows that \ref{subMFCII} of subMFC is fulfilled at $\bar{x}$. Hence, as we suppose that subMFC is violated at $\bar{x}$, there must be a contradiction to \ref{subMFCI} of subMFC for $I= I(x^k,\delta^k)$. Thus, there exists some  $u \geq 0, u \neq 0$ with
\[ 
0 \displaystyle \in \sum_{j\in {I}} u_j \partial g_j(\bar{x}). 
\]
As $I=I(x^k,\delta^k) \subseteq I(\bar{x})$ by the continuity of $g$ and due to $\delta^k \to 0$, it follows for
\[
v_j \coloneqq 
\begin{cases}
u_j & \forall j \in I,  \\
0 & \forall j \in I(\bar{x}) \setminus I
\end{cases}
\]
that the inclusion in \ref{item:quasinorm1} of Definition \ref{def:pq_normality}
holds. Thus, for all $j$ with $v_j>0$ we have $j \in I = I(x^k,\delta^k)$, and hence,
\[ 
v_j g_j(x^k) = v_j \delta^k_j > 0 ,
\]
due to $\delta^k_j > 0$ for all $j=1,\dotsc,q$. Hence, as \ref{item:quasinorm2} of Definition \ref{def:pq_normality} holds as well, quasinormality is violated at $\bar{x}$. By this contradiction, it follows that quasinormality implies subMFC at $\bar{x}$. Finally, Example~\ref{ex:mehrelementiges_subdiff} below shows that this implication is strict.
\end{proof}

In the subsequent example it is shown that subMFC can serve as a condition to identify a KKT point in a situation where pseudonormality and quasinormality fail. 

\begin{example}
\label{ex:mehrelementiges_subdiff}
Consider the Lipschitzian optimization problem
\[
\min \, (x_1+2)^2 + (x_2-1)^2 \quad \text{ s.t. } \quad g_1(x_1,x_2) \coloneqq x_2^3 - 3x_2 - \widetilde{g}_1(x_1) \leq 0,
\]
where
\[ 
\widetilde{g}_1(x_1) \coloneqq \begin{cases}
    x_1^3-3x_1, & \text{for } x_1\in(-\infty, -2) \cup (1,\infty), \\
    -2, & \text{for } x_1\in[-2,1].
  \end{cases} 
\]
The function $g_1$ is not differentiable but locally Lipschitz continuous. Its Clarke subdifferential is given by
\[ 
\partial g_1(x_1,x_2) = 
 \begin{pmatrix} \partial(- \widetilde{g}_1)(x_1)  \\ 3x_2^2-3 \end{pmatrix},\,\,\,
 \partial(- \widetilde{g}_1)(x_1) =
 \begin{cases}
   -3x_1^2+3, & \text{for } x_1\in(-\infty, -2) \cup (1,\infty),  \vspace{2pt} \\
   [-9,0], & \text{for } x_1 = -2, \vspace{2pt} \\
   0, & \text{for } x_1\in(-2,1].
  \end{cases}
\]
It is straightforward to see that the global solution to the problem is $\bar{x} = (\bar{x}_1, \bar{x}_2) = (-2,1)$  
and that it is a KKT point with multiplier $u=0$ and $I(\bar{x}) = \{ 1 \}$. 
We now take a look at the different conditions pseudonormality, quasinormality and subMFC:

\begin{enumerate}
\item First, we consider quasinormality at $\bar{x}$. It follows from $I(\bar{x}) = \{ 1 \}$ and 
\[ 
\partial g_1(\bar{x}_1,\bar{x}_2) = \begin{pmatrix} [-9,0] \\ 0 \end{pmatrix}
\]
that part \ref{item:quasinorm1} of Definition \ref{def:pq_normality} is fulfilled for any $v_1 > 0$. By taking the sequence $\{ (x_1^k,x_2^k) \}_{k=1}^\infty$ with $(x_1^k,x_2^k) \coloneqq \left(-2+\frac{1}{k}, 1+\frac{1}{k}\right) \to (\bar{x}_1,\bar{x}_2)$ as $k\to\infty$, we see that
\begin{align*}
v_1 g_1(x_1^k,x_2^k) &= v_1( (x_2^k)^3 -3 x_2^k +2) 
= v_1\left( \frac{1}{k^3} + \frac{3}{k^2} \right) > 0 \qquad \forall v_1 > 0.
\end{align*} 
Thus, the point $\bar{x}$ does not satisfy quasinormality which implies that pseudonormality can neither be fulfilled.
\item Next, we show that subMFC is fulfilled at $\bar{x}$. The AKKT conditions in $\bar{x}$ for this problem require $(x^k, \delta^k,\epsilon^k) \to (\bar{x},0,0)$ for $k\to\infty$ and, for all $k\in\mathbb{N}$,
\begin{gather*}
\begin{pmatrix} \epsilon^k_1 \\ \epsilon^k_2 \end{pmatrix} \in \begin{pmatrix} 2(x_1^k+2) \\ 2(x_2^k-1) \end{pmatrix} + u^k  \begin{pmatrix} \partial (-\widetilde{g}_1)(x_1^k)  \\ 3(x_2^k)^2-3 \end{pmatrix}, \\
u^k \geq 0,\quad (x_2^k)^3 - 3x_2^k - \widetilde{g}_1(x_1^k)  \leq \delta^k,\quad u^k(x_2^3 - 3x_2 - \widetilde{g}_1(x_1^k)-\delta^k) = 0.
\end{gather*}
It can be easily verified that the sequence
\[ 
\{({x}^k,u^k, \delta^k, \epsilon^k)\}_{k=1}^\infty \coloneqq \left\{\left(\left(-2+\frac{1}{k},1+\frac{1}{k}\right),0,\frac{6}{k^2},\left(\frac{2}{k},\frac{2}{k}\right)\right)\right\}_{k=1}^\infty 
\]
fulfills the above AKKT conditions at $\bar{x}$.
For this sequence, $I(x^k,\delta^k) = \emptyset$ holds for all~$k\in\mathbb{N}$. Thus, we immediately see that subMFC is satisfied at $\bar{x}$ with $I \coloneqq \emptyset$.
\end{enumerate}
\end{example}

\subsection{Local error bound condition}
\label{subsec:errorbound}
It is known that CQs may imply local error bounds to the feasible region if the constraint functions are sufficiently smooth, see, for instance, \cite{RA2} for CPLD, \cite{LM1} for CRCQ and \cite{LM2} for quasinormality. Some of these results can, under additional assumptions, be extended to the nonsmooth case. For example, if the constraint functions are subdifferentially regular at the point of interest $\bar{x}$, which means that $\widehat{\partial}  g_j(\bar{x}) = \partial^L g_j(\bar{x})$ holds for all $j=1,\dotsc,q$, the local error bound condition is implied by CPLD \cite[Theorem 3.2]{MX1} and by quasinormality \cite[Theorem 5]{JY11}. In \cite[Corollary 5.3]{LG1}, it was shown that for the latter result the additional assumption on subdifferential regularity of the constraint functions can even be dropped. Thus, as all of the standard CQs considered so far are stronger than subMFC, we want to take a look at the local error bound condition now. To this end, recall that the feasible set of \eqref{P} is denoted by $D$ and that $d_D(x)$ stands for the distance of $x\in\mathbb{R}^n$ to the set $D$. 

\begin{definition}
The \textit{local error bound condition} is said to hold at $\bar{x}$ feasible for \eqref{P} if there exist
constants $c>0$ and $r>0$ such that 
\[ 
d_D(x) \leq c \norm{g_+(x)}  \quad\forall x\in B(\bar{x},r). 
\]
\end{definition}

It turns out that subMFC can hold even if the local error bound condition does not. On the other hand, the local error bound condition may be satisfied even if subMFC is violated.

\begin{proposition}
Let $\bar{x}$ be feasible for \eqref{P}. Then, subMFC at $\bar{x}$ and the local error bound condition at $\bar{x}$ are independent conditions.
\end{proposition}
\begin{proof}
See Appendix \ref{sec:App_subMFC_errorbound}.
\end{proof}

\subsection{Cone-continuity property}
Next, we want to look at the weakest CQ which implies that an AKKT point is a KKT point. This CQ was established in \cite{RA5} for smooth optimization problems and called CCP. It has been extended to problems with locally Lipschitz continuous functions in \cite{PM1}, where it was renamed as asymptotic regularity. 

\begin{definition}
We say that \textit{CCP} holds at a point $\bar{x}$ feasible for \eqref{P} if the set-valued mapping $K:\mathbb{R}^n \rightrightarrows \mathbb{R}^{n}$ defined by
\[ 
K(x) = \left\{ \sum_{j \in I(\bar{x})} v_j \partial g_j(x) \mid v_j\geq 0 \right\} 
\]
is outer semicontinuous at $\bar{x}$, i.e., if
\begin{alignat*}{2}
\limsup_{x \to \bar{x}} K(x) &\coloneqq \{ \bar{z} \in\mathbb{R}^n \mid \,\,\exists \{ (x^k, z^k) \}_{k=1}^\infty \subset \mathbb{R}^{n+n}:\,\, && (x^k,z^k) \to (\bar{x},\bar{z}), z^k\in K(x^k)\,\,\, \forall k\in\mathbb{N}\}\\
&\subset K(\bar{x}).&& \numberthis\label{eq:CCP} 
\end{alignat*}
\end{definition}

In \cite[Lemma 3.12]{PM1}, it is shown that CCP is implied by CPLD, and it follows from \cite{RA5} that this implication is strict. Moreover, in \cite{RA5}, it is shown that CCP neither implies nor is implied by quasinormality.
In \cite{RA5,EB1}, it is proven that CCP is the weakest possible strict CQ in the smooth case and in \cite{PM1}, a similar claim is made for the nonsmooth case. Due to Theorem \ref{thm:bqc_implies_kkt}, we obtain that subMFC also possesses the most important property of strict CQs, which is guaranteeing that an AKKT point is already a KKT point. At the same time, subMFC neither implies nor is implied by CCP, as we will see in the next proposition. This is reasonable as the key point is to recognize that subMFC is not a CQ, as part \ref{subMFCII} of its definition also depends on the given objective function. In contrast, strict CQs have in common that they can be stated without depending on the specific AKKT point sequence. It turns out that dropping this last property enables subMFC to be such a powerful condition: Indeed, in subMFC, only \textit{one} AKKT point sequence has to be considered, and the defining property~\ref{subMFCI} only needs to be valid for this particular sequence. Due to this link of the AKKT point sequence with an additional property, subMFC can hold even if usual strict CQs fail, while stile maintaining their most important property. The announced proposition illustrates this circumstance.

\begin{proposition}
Let $\bar{x}$ be feasible for \eqref{P}. Then, subMFC at $\bar{x}$ and CCP at $\bar{x}$ are independent conditions.
\end{proposition}
\begin{proof}
See Appendix \ref{sec:App_subMFC_CCP}.
\end{proof}

\subsection{Abadie constraint qualification and Guignard constraint qualification}
\label{sec:acq_gcq}
Now, we consider the Abadie CQ (ACQ), which is not a strict CQ, but comparatively weak as it is implied by both CCP \cite{RA5} and quasinormality \cite{DB2} for the differentiable case. 
Moreover, ACQ is equivalent to the local error bound condition if the constraint functions $g_j$, $j=1,\dots,q$, are differentiable and convex \cite[Theorem 3.5]{WL1}. As similar results do not exist for the nonsmooth case, we thus investigate the relation of subMFC to ACQ separately in this section. 
To this end, we use the definition of the nonsmooth ACQ presented in \cite{R2}. It is based on the \textit{contingent cone} to the feasible set $D$ of problem \eqref{P} at $\bar{x}\in\mathbb{R}^n$, given by
\[ 
T_D(\bar{x}) \coloneqq \left\{ d\in\mathbb{R}^n \mid \exists \{t_k\}_{k=1}^\infty,\{d^k\}_{k=1}^\infty:  t_k \downarrow 0, d^k \to d \text{ with } \bar{x}+t_kd^k\in D\,\,\, \forall k\in\mathbb{N} \right\}, 
\]
and the \textit{cone of locally constrained directions} to $D$ at $\bar{x}\in\mathbb{R}^n$, which is
\[ 
L_D(\bar{x}) \coloneqq \left\{ d\in\mathbb{R}^n \mid s^\top d \leq 0\quad \forall s\in \bigcup_{j\in I(\bar{x})} \partial g_j(\bar{x}) \right\}. 
\]

\begin{definition}
\textit{ACQ} is said to hold at a point $\bar{x}$ feasible for \eqref{P} if
\[ 
L_D(\bar{x}) \subseteq T_D(\bar{x}). 
\]
\end{definition}

Moreover, we consider the Guignard CQ (GCQ), which is even weaker than ACQ. Similar to the smooth case, it is derived by replacing the contingent cone in the definition of ACQ with a suitable enlargement, as presented in \cite{R2}.

\begin{definition}
\textit{GCQ} is said to hold at a point $\bar{x}$ feasible for \eqref{P} if
\[ 
L_D(\bar{x}) \subseteq \cl \conv T_D(\bar{x}),
\]
where $\cl \conv T_D(\bar{x})$ is the closure of the convex hull of $T_D(\bar{x})$.
\end{definition}

Under mild assumptions, it can be shown that these nonsmooth extensions of the smooth ACQ and GCQ are indeed constraint qualifications for \eqref{P}, see \cite[Theorem 2.3]{R2} for details. In the next proposition, it is shown that subMFC is independent of both of these weak CQs.

\begin{proposition}
Let $\bar{x}$ be feasible for \eqref{P}. Then, subMFC at $\bar{x}$ and ACQ/GCQ at~$\bar{x}$ are independent conditions.
\end{proposition}
\begin{proof}
See Appendix \ref{sec:AppA_subMFC_ACQGCQ}. 
\end{proof}

\subsection{Calmness in the sense of Clarke}
\label{subsec:calmness}
Finally, we consider calmness in the sense of Clarke \cite[Definition 6.4.1]{FC1}.

\begin{definition}
Let $\bar{x}$ solve \eqref{P}. Then, \eqref{P} is said to be \textit{calm} at $\bar{x}$ if there exist constants $r>0$ and $\mu>0$ such that, for all $\alpha \in B(0,r)$ and for all $x \in B(\bar{x},r)$ with
\[ 
g(x) + \alpha \leq 0, 
\] 
it follows that 
\[ 
f(x) - f(\bar{x}) + \mu \norm{\alpha} \geq 0. 
\]
\end{definition}

Note that calmness can also be defined for $\bar{x}$ being a local minimizer of \eqref{P}. It is known that calmness is weaker than MFCQ at local minimizers of \eqref{P}, see \cite[Corollary 5 of Theorem 6.5.2]{FC1}. Like subMFC, calmness is a problem qualification and not a CQ, as it involves the objective function. This allows calmness to be comparatively weak, but even this condition can be violated while subMFC still holds, as the following proposition shows.

\begin{proposition}
Let $\bar{x}$ be feasible for \eqref{P}. Then, subMFC at $\bar{x}$ does not imply calmness in the sense of Clarke at~$\bar{x}$.
\end{proposition}
\begin{proof}
See Appendix \ref{sec:AppC_subMFC_calmness}.
\end{proof}

\subsection{Summary of the comparisons}
\label{subsec:summarygraph}
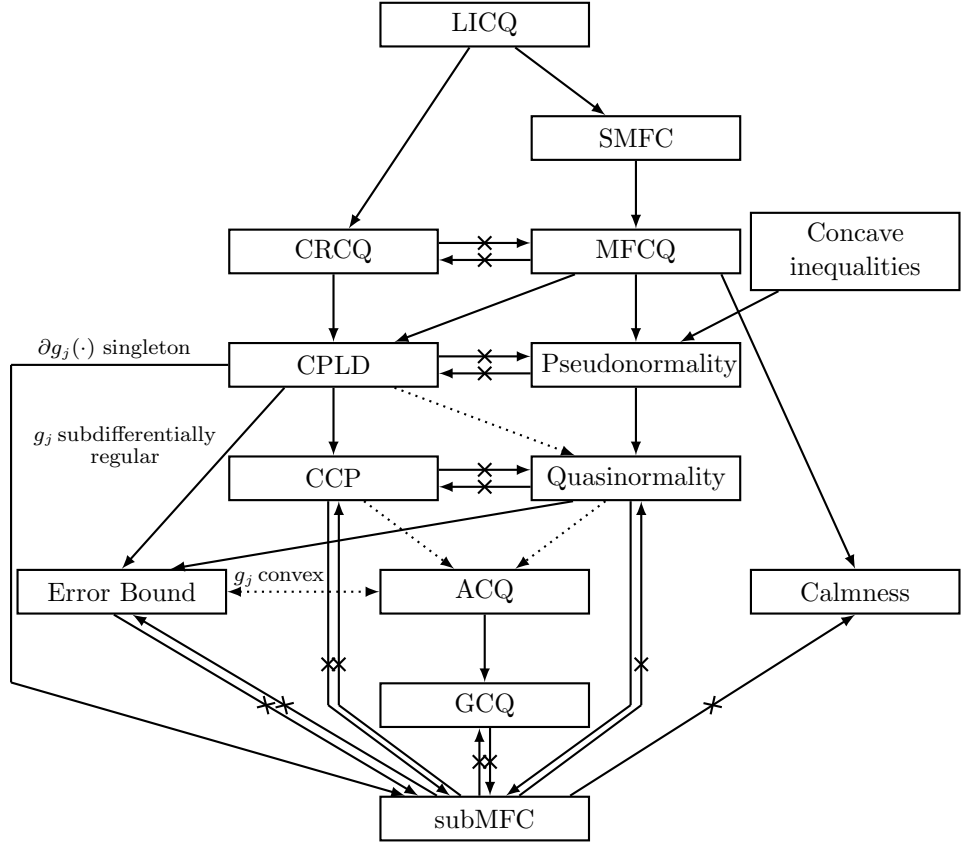
\begin{figure}[H]
\begin{center}
\begin{tikzpicture}
\tikzset{
 block/.style ={rectangle, draw=black, thick, text width=6.5em,align=center, minimum height=1.5em},
 eblock/.style ={align=center, minimum height=1.5em},
 arrow/.style = {thick,-{latex[scale=4]}},
 smootharrow/.style = {thick,dotted,-{latex[scale=4]}},
 smoothequi/.style = {thick,dotted,{latex[scale=4]}-{latex[scale=4]}},
 strike through/.style={postaction=decorate, decoration={markings,mark=at position 0.5 with {\draw[-] (-2.5pt,-2.5pt)--(2.5pt,2.5pt); \draw[-,rotate=90] (-2.5pt,-2.5pt)--(2.5pt,2.5pt);}}},
 strike throughtwo/.style={arrow, postaction=decorate, decoration={markings,mark=at position 0.2 with {\draw[-] (-2.5pt,-2.5pt)--(2.5pt,2.5pt); \draw[-,rotate=90] (-2.5pt,-2.5pt)--(2.5pt,2.5pt);}}},
 strike throughthree/.style={thick,postaction=decorate, decoration={markings,mark=at position 0.8 with {\draw[-] (-2.5pt,-2.5pt)--(2.5pt,2.5pt); \draw[-,rotate=90] (-2.5pt,-2.5pt)--(2.5pt,2.5pt);}}}
}
\draw (0, 1.5) node[block] (A) {\small LICQ};
\draw (2, 0) node[block] (AB) {\small SMFC};
\draw (-2,-1.5) node[block] (B) {\small CRCQ};
\draw (2,-1.5) node[block] (C) {\small MFCQ};
\draw (4.9,-1.5) node[block] (D) {\small Concave \\ inequalities};
\draw (-2,-3) node[block] (E) {\small CPLD};
\draw (2,-3) node[block] (F) {\small Pseudonormality};
\draw (2,-4.5) node[block] (G) {\small Quasinormality};
\draw (0,-9) node[block] (H) {\small subMFC};
\draw (0,-5) node (I) {};
\draw (-4.8,-6) node[block] (J) {\small Error Bound};
\draw (4.9,-6) node[block] (K) {\small Calmness};
\draw (-2,-4.5) node[block] (L) {\small CCP};
\draw (0,-6) node[block] (M) {\small ACQ};
\draw (-6.4,-3) node[eblock] (N) {};
\draw (-6.4,-7.2) node[eblock] (O) {};
\draw (0,-7.5) node[block] (P) {\small GCQ};
\draw [arrow] (A) -- (B);
\draw [arrow] (A) -- (AB);
\draw [arrow] (AB) -- (C);
\draw [->, arrow, strike through] ([yshift= 3.2pt] B.east) -- ([yshift= 3.2pt] C.west);
\draw [<-, arrow, strike through] ([yshift= -3.2pt] C.west) -- ([yshift= -3.2pt] B.east);
\draw [arrow] (B) -- (E); 
\draw [arrow] (C) -- (E);
\draw [arrow] (C) -- (F);
\draw [arrow] (D) -- (F); 
\draw [->, arrow, strike through] ([yshift= 3.2pt] E.east) -- ([yshift= 3.2pt] F.west);
\draw [<-, arrow, strike through] ([yshift= -3.2pt] F.west) -- ([yshift= -3.2pt] E.east);
\draw [smootharrow] (E) -- (G);
\draw [arrow] (E.205) -- node[above left, xshift=9pt] {$\substack{g_j\, \text{subdifferentially} \\ \text{regular}}$} (J.north);
\draw [arrow] (F) -- (G);
\draw [thick] (E.west) -- node[above, xshift=-2pt, yshift=-2pt] {$\substack{ \\ \partial g_j(\cdot)\,\, \text{singleton}}$} (N.east);
\draw [thick] (N.east) -- (O.east);
\draw [arrow] (O.east) -- (H); 
\draw [thick] ([xshift= -2pt] G.south) -- ([xshift= -2pt] 2,-7.5);
\draw [arrow] ([xshift= -2pt] 2,-7.5) -- ([xshift= 8pt] H.north);
\draw [thick] ([xshift= 13pt] H.north) -- ([xshift= 2pt] 2,-7.5);
\draw [<-, arrow, strike throughtwo] ([xshift= 2pt] 2,-7.5) -- ([xshift= 2pt] G.south);
\draw [arrow] (C.345) -- (K.north); 
\draw [->, arrow, strike through] (H.15) -- (K.south);
\draw [arrow] (E) -- (L); 
\draw [thick] ([xshift= -9pt] H.north) -- ([xshift= 2pt] -2,-7.5);
\draw [strike throughtwo] ([xshift= 2pt] -2,-7.5) -- ([xshift= 2pt] L.south);
\draw [strike throughthree] ([xshift= -2pt] L.south) -- ([xshift= -2pt] -2,-7.5);
\draw [thick, arrow] ([xshift= -2pt] -2,-7.5) -- ([xshift= -13pt] H.north);
\draw [arrow, strike through] ([xshift= -18pt] H.north) -- ([xshift= 4pt] J.south);
\draw [arrow, strike through] ([xshift= -3pt] J.south) -- ([xshift= -25pt] H.north);
\draw [arrow] (M) -- (P); 
\draw [arrow, strike through] ([xshift= -2pt] H.north) -- ([xshift= -2pt] P.south);
\draw [arrow, strike through] ([xshift= 2pt] P.south) -- ([xshift= 2pt] H.north);
\draw [smootharrow] (L) -- (M);
\draw [smootharrow] (G) -- (M); 
\draw [->, arrow, strike through] ([yshift= -3.2pt] G.west) -- ([yshift= -3.2pt] L.east);
\draw [->, arrow, strike through] ([yshift= 3.2pt] L.east) -- ([yshift= 3.2pt] G.west);
\draw [arrow] (G.200) -- (J.25);
\draw [smoothequi] (J) -- node[above, xshift=-9pt, yshift=-2pt] {$\substack{\\ g_j\, \text{convex}}$} (M);
\end{tikzpicture}
\end{center}
\caption{Strict implications between different conditions considered at a local minimizer of the inequality constrained problem \eqref{P} are represented by solid arrows. Additional strict implications and one equivalence in case of differentiability of the constraint functions are visualized by dotted arrows. A crossed out arrow symbolizes that the condition at the start of the arrow does not imply the one at its end.}
\label{fig:CQ_implications}
\end{figure}

In Figure \ref{fig:CQ_implications}, the implications between the conditions considered in this paper are summarized. In particular, we want to emphasize that all implications are strict. All depicted relations can be verified using the results presented throughout Section \ref{sec:other_conds}, together with the examples and results listed for completeness in Appendix \ref{sec:examples_strictness}.

\section{Application to bilevel optimization problems}
In this section, we tailor subMFC to the lower-level value function reformulation
\begin{equation}
\min_{x,y}\,\, F(x,y) \quad\text{s.t.}\quad G(x,y)\leq 0,\quad \hat{g}(x,y)\leq 0,\quad \hat{f}(x,y) - \varphi(x) \leq 0 \tag{LLVFR}\label{LLVFR}
\end{equation}
of the optimistic bilevel optimization problem \eqref{BP}. We want to emphasize that we do not require any additional assumptions besides the local Lipschitz continuity of the bilevel problem data functions $F$, $G$, $\hat{f}$ and $\hat{g}$ and the well-definedness and local Lipschitz continuity of the lower-level value function $\varphi$, which are rather weak assumptions in bilevel optimization. To guarantee the local Lipschitz continuity of the lower-level value function $\varphi$ near some $x$, many different sufficient conditions are known, see, for example, \cite[Proposition 1.2]{JY6}, \cite[Corollary 1 of Theorem 6.5.2]{FC1} or \cite[Theorem 3.4]{SD1}, and the references mentioned therein.

For a feasible point $(x,y)$ of problem \eqref{LLVFR}, we denote the index sets of active constraints by
\begin{align*}
I^G(x,y) &\coloneqq \{i \in \{1,\dotsc,p\} \mid G_i(x,y)=0\}, \\
I^{\hat{g}}(x,y) &\coloneqq \{j \in \{1,\dotsc,r\} \mid \hat{g}_j(x,y)=0\}, \\
I^{\hat{f}-\varphi}(x,y) &\coloneqq \{ 1 \mid \hat{f}(x,y) - \varphi(x) = 0 \} = \{1 \}.
\end{align*}
Note that the last equation is valid because, by definition of $\varphi$, $\hat{f}-\varphi$ is inherently active in any point that is feasible for \eqref{LLVFR}. The latter fact also leads to standard CQs failing for every feasible point of \eqref{LLVFR}, which is, among others, one reason to relax problem \eqref{LLVFR} by
\begin{equation}
\min_{x,y}\,\, F(x,y) \quad\text{s.t.}\quad G(x,y)\leq 0,\quad \hat{g}(x,y)\leq 0,\quad \hat{f}(x,y) -\varphi(x) \leq \gamma \label{LLVFe}\tag{\mbox{LLVFR$_\gamma$}} 
\end{equation}
with a perturbation $\gamma>0$ on the lower-level value function constraint \cite{LL1,GL1,JY4}.
In contrast to \eqref{LLVFR}, the constraint including the lower-level value function may or may not be active in any point that is feasible for \eqref{LLVFe}, and hence, we define
\[ 
I^{\hat{f}-\varphi}({x},{y},\gamma) \coloneqq \begin{cases} \{ 1 \}, &  \text{if } \hat{f}(x,y) -\varphi(x) - \gamma = 0, \\ \emptyset, & \text{if } \hat{f}(x,y) -\varphi(x) - \gamma < 0. \end{cases}
\]

For several subclasses of the problem class \eqref{LLVFe}, standard CQs can hold at feasible points for a fixed value $\gamma>0$ \cite{LL1,GL1,JY4}. Thus, our aim is to build a connection between KKT points of \eqref{LLVFe} and meaningful points of \eqref{LLVFR}, for which we have the following first result.

\begin{proposition}
\label{prop:conv_fj}
Let $\{ \gamma_k \}_{k=1}^\infty$ be a sequence of positive reals such that $\gamma_k \downarrow 0$ as \mbox{$k\to\infty$}. 
Assume that the point $(x^{k},y^{k})$ is a KKT point of \hyperref[LLVFe]{\textup{(LLVFR$_{\gamma_k}$)}} for each $k\in\mathbb{N}$. 
Then, any accumulation point of $\{(x^{k},y^{k})\}_{k=1}^\infty$ is feasible for \eqref{LLVFR}.
\end{proposition}

This result can directly be verified using the continuity of the functions $G$, $\hat{f}$, $\hat{g}$ and~$\varphi$. Knowing that a sequence of KKT points of \eqref{LLVFe} converges to a feasible point of \eqref{LLVFR} for $\gamma \downarrow 0$, we want to investigate next under which conditions more than just feasibility can be transferred to the limit point, i.e., when is this point an FJ point, an AKKT point or a KKT point for \eqref{LLVFR}? For the first part, we obtain the following answer which aligns with \cite[Proposition 3.2]{JY3}. To keep the paper self-contained, we also present a short proof.

\begin{proposition}
\label{prop:llvf_fj}
Any point that is feasible for \eqref{LLVFR} is an FJ point of \eqref{LLVFR} with $a=0$.
\end{proposition}

\begin{proof}
By definition of $\varphi$ we have $\hat{f}(x,y)\geq \varphi(x)$ for all $(x,y)\in\{(x,y)\mid \hat{g}(x,y)\leq 0\}$ and consequently, $\hat{f}(x,y) = \varphi(x)$ holds at all points feasible for \eqref{LLVFR}. Thus, any point $(\bar{x},\bar{y})$ feasible for \eqref{LLVFR} is a global solution to the problem
\[ 
\min_{x,y} \hat{f}(x,y)-\varphi(x) \quad\text{s.t.}\quad G(x,y)\leq 0,\quad \hat{g}(x,y)\leq 0. 
\]
As $\varphi$ is Lipschitz-continuous, we can apply the necessary FJ optimality conditions from Definition \ref{def:fj} to the above problem and obtain, for any $(\bar{x},\bar{y})$ feasible for \eqref{LLVFR}, the existence of multipliers $w\geq0$, $u\geq 0$, $v\geq 0$ with $(w,u,v) \neq (0,0,0)$ such that 
\[ 
0 \in w \partial (\hat{f}-\varphi)(\bar{x},\bar{y}) + \sum_{i\in I^G(\bar{x},\bar{y})} u_i \partial G_i(\bar{x},\bar{y}) + \sum_{j\in I^g(\bar{x},\bar{y})} v_j \partial \hat{g}_j(\bar{x},\bar{y}). 
\]
Consequently, the relation
\[ 
0 \in a \partial F(\bar{x},\bar{y})+ \sum_{i\in I^G(\bar{x},\bar{y})} u_i \partial G_i(\bar{x},\bar{y}) + \sum_{j\in I^g(\bar{x},\bar{y})} v_j \partial \hat{g}_j(\bar{x},\bar{y}) + w \partial (\hat{f}-\varphi)(\bar{x},\bar{y}) 
\]
holds for $a =0$, which shows that $(\bar{x},\bar{y})$ is always an FJ point of \eqref{LLVFR}, but not necessarily a KKT point. 
\end{proof}
  
Unfortunately, the combination of Proposition \ref{prop:conv_fj} and Proposition \ref{prop:llvf_fj} shows that, without further assumptions, all accumulation points of sequences of KKT points of \eqref{LLVFe}, computed as $\gamma\downarrow 0$, are only FJ points of \eqref{LLVFR}. 
As any feasible point shares the same property, we aim to ensure the convergence to a more meaningful point. Thus, naturally the question arises whether being an AKKT point is a more restrictive property than being an FJ point. The next example shows that this is indeed the case.

\begin{example}
\label{ex:strict_implications}
Consider the optimistic bilevel optimization problem
\begin{alignat*}{2}
\min_{x,y} -y\quad \text{ s.t. }\quad y\in S(x) \coloneqq \argmin_y \{xy \mid\,\, && -y \leq 0, 
\,\, y-1 \leq 0 \}. 
\end{alignat*}
The lower-level value function $\varphi$ is given by
\[ 
\varphi(x) = \begin{cases} x, & \text{if } x \leq 0, \\ 0, & \text{if } x > 0. \end{cases}
\]
This yields the subdifferential 
\vspace{-0.2em}
\[
\partial (\hat{f}-\varphi)(x,y)=  \begin{pmatrix} h(x,y) \\ x \end{pmatrix}
\quad\text{with}\quad h(x,y)  = \begin{cases}
    y-1, & \text{for } x<0, \\
    [y-1,y], & \text{for } x = 0, \\
    y, & \text{for } x>0
  \end{cases}
\]
and the lower-level value function reformulation 
\[ 
\min_{x,y}\,\, -y \quad\text{s.t.}\quad -y\leq 0,\quad y-1 \leq 0,\quad xy - \varphi(x) \leq 0.
\]
Clearly, the point $(\bar{x},\bar{y}) \coloneqq \left(0,\frac{1}{2}\right)$ is feasible to the problem. 
The AKKT conditions \eqref{eq:akkt1}--\eqref{eq:akkt5} in $(\bar{x},\bar{y})$ for the reformulation require $(x^k, y^k, \delta^k, \epsilon^k) \to (\bar{x},\bar{y}, 0, 0)$ as $k\to\infty$ and
\begin{gather*}
\epsilon^k \in \begin{pmatrix} 0 \\ -1 \end{pmatrix} + v^k_1  \begin{pmatrix} 0 \\ -1 \end{pmatrix} + v^k_2 \begin{pmatrix} 0 \\ 1 \end{pmatrix} + w^k  \begin{pmatrix} h(x^k,y^k) \\ x^k \end{pmatrix}, \numberthis \label{eq:ex_fj_akkt} \\
\begin{alignat*}{4}
&& v^k_1 \geq 0,&& -y^k \leq \delta^k_1,&& v^k_1 (-y^k-\delta^k_1) = 0, \phantom{aaaaaaaa,}\\
&& v^k_2 \geq 0,&&\quad y^k-1 \leq \delta^k_2,&&\quad v^k_2 (y^k-1 -\delta^k_2) = 0, \phantom{aaaaaaaa,} \\
&& w^k \geq 0,&&\quad x^ky^k-\varphi(x^k) \leq \delta^k_3,&&\quad w^k (x^ky^k-\varphi(x^k) - \delta^k_3 ) = 0 \phantom{,aaaaaaaa,}
\end{alignat*}
\end{gather*}
for all $k\in\mathbb{N}$.
Due to $y^k \to \frac{1}{2}$ and $\delta^k \to 0$, we immediately obtain $v^k_1 = v^k_2 = 0$ for large $k\in\mathbb{N}$. Thus, the second component in \eqref{eq:ex_fj_akkt} reads as $\epsilon^k_2 = -1 + w^k x^k$. As $x^k \to 0$ and $\epsilon^k \to 0$ are required, it follows that this equation cannot be fulfilled. Hence, the point $(\bar{x},\bar{y})$ is feasible (and thus an FJ point) but not an AKKT point.\\
Moreover, the sequence
\vspace{-0.2em}
\[ 
\{ (x^k,y^k,v^k,w^k,\delta^k,\epsilon^k) \}_{k=1}^\infty \coloneqq \left\{ \left( \frac{1}{k},0,\left(0,0\right),k,(0,0,0),(0,0) \right) \right\}_{k=1}^\infty 
\]
can be used to show that the point $(\widetilde{x},\widetilde{y}) \coloneqq (0,0)$ is an AKKT point, as it fulfills the above system. However, it can be checked that $(\widetilde{x},\widetilde{y})$ is not a KKT point, which is important for the next result.
\end{example}

With Example \ref{ex:strict_implications}, we can establish the following strict relations between the different stationarity concepts considered in this paper.

\begin{proposition}
\label{prop:strict_implications}
Let $(\bar{x},\bar{y})$ be feasible for \eqref{LLVFR}.
Then, the following implications are strict:
\[ (\bar{x},\bar{y}) \text{ is a KKT point} \Longrightarrow (\bar{x},\bar{y}) \text{ is an AKKT point} \Longrightarrow (\bar{x},\bar{y}) \text{ is an FJ point}. \]
\end{proposition} 
\begin{proof}
From Definitions \ref{def:kkt} and \ref{def:akkt} it follows that, even for the general problem \eqref{P}, any KKT point $\bar{x}$ with multiplier $u$ is an AKKT point with the corresponding sequence $(x^k,u^k,\delta^k,\epsilon^k)\coloneqq (\bar{x},u,0,0)$. The second implication holds for \eqref{LLVFR} due to Proposition \ref{prop:llvf_fj} and the fact that any AKKT point is in particular feasible, as the constraint functions are continuous. The strictness of the implications then follows from Example \ref{ex:strict_implications}. 
\end{proof}

As the following proposition shows, it turns out that no additional assumptions are required to guarantee that an accumulation point of a sequence of KKT points of problem \eqref{LLVFe} for $\gamma\downarrow 0$ fulfills the stronger requirement of being an AKKT point of \eqref{LLVFR}.

\begin{proposition}
\label{prop:conv_akkt}
Let $\{ \gamma_k \}_{k=1}^\infty$ be a sequence of positive reals such that $\gamma_k \downarrow 0$ as $k\to\infty$. 
Assume that the point $(x^{k},y^{k})$ is a KKT point of \hyperref[LLVFe]{\textup{(LLVFR$_{\gamma_k}$)}} for each $k\in\mathbb{N}$. 
Then, any accumulation point of $\{(x^{k},y^{k})\}_{k=1}^\infty$ is an AKKT point of \eqref{LLVFR}.
\end{proposition}

\begin{proof}
Consider an accumulation point $(\bar{x},\bar{y})$ of the sequence $\{(x^{k},y^{k})\}_{k=1}^\infty$. Without loss of generality, we may assume that the whole sequence $\{(x^{k},y^{k})\}_{k=1}^\infty$ converges to $(\bar{x},\bar{y})$. As $(x^{k},y^{k})$ is a KKT point of \hyperref[LLVFe]{(LLVFR$_{\gamma_k}$)}, there are multipliers $u^k$, $v^k$ and $w^k$ such that
\begin{gather*}
0 \in \partial F(x^{k},y^{k}) + \partial G(x^{k},y^{k})u^k + \partial \hat{g}(x^{k},y^{k}) v^k + w^k \partial(\hat{f}-\varphi)(x^{k},y^{k}), \\
\begin{alignat*}{4}
&& u^k \geq 0,&& G(x^{k},y^{k}) \leq 0,\,\,\,&& (u^k)^\top G(x^{k},y^{k}) = 0, \phantom{aiaaaaa} \\
&& v^k \geq 0,&& \hat{g}(x^{k},y^{k}) \leq 0,\,\,\,&& (v^k)^\top \hat{g}(x^{k},y^{k}) = 0, \phantom{aiaaaaa} \\
&& w^k \geq 0,&&\quad \hat{f}(x^{k},y^{k})-\varphi({x}^k) \leq \gamma_k,&&\quad w^k (\hat{f}(x^{k},y^{k})-\varphi(x^{k})-\gamma_k) = 0 \phantom{,}\phantom{aiaaaaa}
\end{alignat*}
\end{gather*}
holds for all $k\in\mathbb{N}$. By the assumptions, we further have
\begin{equation*}
(x^k, y^k, \gamma_k) \to (\bar{x},\bar{y},0) \quad \text{as}\,\,\, k\to\infty.
\end{equation*}
Consequently, the definition of an AKKT point is fulfilled at $(\bar{x},\bar{y})$ with the sequence $\{(\widehat{x}^k,\widehat{u}^k,\widehat{\delta}^k,\widehat{\epsilon}^k)\}_{k=1}^\infty$ for $\widehat{x}^k \coloneqq (x^k,y^k)$, $\widehat{u}^k  \coloneqq (u^k,v^k,w^k)$, $\widehat{\delta}^k \coloneqq (0,0,\gamma_k)$ and $\widehat{\epsilon}^k \coloneqq (0,0)$. 
\end{proof}

Finally, we are interested in finding a suitable condition under which we obtain a KKT point of \eqref{LLVFR} when looking at a sequence of KKT points of problem \eqref{LLVFe} for $\gamma\downarrow 0$. 
As a consequence of Proposition \ref{prop:conv_akkt}, we can use the results from Section \ref{sec:subMFC} to establish a connection between KKT points of \eqref{LLVFe} and KKT points of \eqref{LLVFR}. To this end, we introduce a new problem qualification for the problem class \eqref{LLVFR} by tailoring Definition~\ref{def:subMFC}.

\begin{definition}
\label{def:subMFC_bilevel}
Let $(\bar{x},\bar{y})$ be feasible for \eqref{LLVFR}. We say that the \textit{LLVFR-subMFC} holds at $(\bar{x},\bar{y})$ if there exist $I_1 \subseteq I^G(\bar{x},\bar{y})$, $I_2 \subseteq I^{\hat{g}}(\bar{x},\bar{y})$ and $I_3 \subseteq I^{\hat{f}-\varphi}(\bar{x},\bar{y}) = \{ 1 \}$ such that the following two conditions are satisfied:
\begin{enumerate}[label=(\roman*)]
\item \label{subMFCI_bilevel} Either $I_1 = I_2 = I_3 = \emptyset$, or it holds for all $u\geq 0$, $v\geq 0$ and $w \geq 0$ with $(u,v,w) \neq (0,0,0)$ that
\vspace{-0.6em}
\[ 
0 \notin \sum_{i\in I_1} u_i \partial G_i(\bar{x},\bar{y}) + \sum_{j\in I_2} v_j \partial \hat{g}_j(\bar{x},\bar{y}) + \begin{cases} w \partial(\hat{f}-\varphi)(\bar{x},\bar{y}), & \text{if } I_3 = \{1\}, \\ 0, & \text{if } I_3 = \emptyset. \end{cases}
\vspace{-0.2em} 
\]
\item \label{subMFCII_bilevel} There exists a sequence $\{(x^k, y^k, u^k, v^k, w^k, \delta^k, \hat{\delta}^k, \gamma_k, \epsilon^k)\}_{k=1}^\infty$ that fulfills 
\begin{gather}
\epsilon^k \in \partial F(x^{k},y^{k}) + \partial G(x^{k},y^{k})u^k + \partial \hat{g}(x^{k},y^{k})v^k + w^k  \partial(\hat{f}-\varphi)(x^{k},y^{k}), \vspace{-1.35em} \label{eq:akkt1_bilevel} \\
\begin{alignat}{6}
&& u^k \geq 0,&& G(x^{k},y^{k}) &\leq \delta^k,\,\,\,& (u^k)^\top (G(x^{k},y^{k})-\delta^k) = 0,\hspace{0.3em} \numberthis \\
&& v^k \geq 0,&& \hat{g}(x^{k},y^{k}) &\leq \hat{\delta}^k,\,\,\,& (v^k)^\top (\hat{g}(x^{k},y^{k}) - \hat{\delta}^k) = 0, \hspace{0.3em} \numberthis \\
&& w^k \geq 0,&&\quad \hat{f}(x^{k},y^{k})-\varphi({x}^k) &\leq \gamma_k,&\quad w^k (\hat{f}(x^{k},y^{k})-\varphi(x^{k})-\gamma_k) = 0 \phantom{,} \hspace{0.3em} \numberthis
\end{alignat}
\end{gather}
for all $k\in\mathbb{N}$ and
\vspace{-0.4em}
\begin{equation}
(x^k, y^k, \delta^k, \hat{\delta}^k, \gamma_k, \epsilon^k) \to (\bar{x},\bar{y}, 0, 0, 0, 0) \quad \text{as}\,\,\, k\to\infty \label{eq:akkt5_bilevel}
\vspace{-0.3em}
\end{equation} 
with
\vspace{-0.3em}
\begin{align*}
I_1 &= I^{G}(x^k,y^k,\delta^k) \coloneqq \{ i\in \{1, \dots, p\} \mid G_i(x^k,y^k) = \delta_j^k \},\\
I_2 &= I^{\hat{g}}(x^k,y^k,\hat{\delta}^k) \coloneqq \{ j\in \{1, \dots, r\} \mid \hat{g}_j(x^k,y^k) = \hat{\delta}_j^k \}, \\
I_3 &= I^{\hat{f}-\varphi}(x^k,y^k,\gamma_k).
\end{align*}
\end{enumerate}
\end{definition}

In analogy to Remark \ref{rem:seq_const}, we assume without loss of generality that the sets $I^{G}(x^k,y^k,\delta^k)$, $I^{\hat{g}}(x^k,y^k,\hat{\delta}^k)$ and $I^{\hat{f}-\varphi}(x^k,y^k,\gamma_k)$ occurring in Definition~\ref{def:subMFC_bilevel}\,\ref{subMFCII_bilevel} are constant for all $k\in\mathbb{N}$, which is reasonable as it is always possible to consider a suitable infinite subsequence. Similarly to Definition \ref{def:subMFC}\,\ref{subMFCII}, Definition~\ref{def:subMFC_bilevel}\,\ref{subMFCII_bilevel} requires that $(\bar{x},\bar{y})$ is an AKKT point of \eqref{LLVFR} whereas part \ref{subMFCI_bilevel} can be interpreted as a weakened MFCQ at $(\bar{x},\bar{y})$. If $(\bar{x},\bar{y})$ is a local minimizer, there always exists at least one sequence that is a suitable candidate for Definition \ref{def:subMFC_bilevel}\,\ref{subMFCII_bilevel}, as a local minimizer of \eqref{LLVFR} is always an AKKT point by Theorem~\ref{thm:locmin_akkt}.
The next theorem is a direct consequence of Theorem \ref{thm:bqc_implies_kkt}.

\begin{theorem}
\label{thm:llvfrsubMFC_implies_kkt}
If $(\bar{x},\bar{y})$ is feasible for \eqref{LLVFR} and LLVFR-subMFC holds at $(\bar{x},\bar{y})$, then $(\bar{x},\bar{y})$ is a KKT point of \eqref{LLVFR}.
\end{theorem}

Finally, combining Proposition \ref{prop:conv_fj} and Theorem \ref{thm:llvfrsubMFC_implies_kkt} directly yields the desired result.

\begin{corollary}
Let $\{ \gamma_k \}_{k=1}^\infty$ be a sequence of positive reals such that $\gamma_k \downarrow 0$ as $k\to\infty$. 
Assume that the point $(x^{k},y^{k})$ is a KKT point of \hyperref[LLVFe]{\textup{(LLVFR$_{\gamma_k}$)}} for each $k\in\mathbb{N}$. 
Then, any accumulation point of $\{(x^{k},y^{k})\}_{k=1}^\infty$ that fulfills LLVFR-subMFC is a KKT point of \eqref{LLVFR}.
\end{corollary}

\subsection{Comparison to partial calmness}
Partial calmness was introduced in \cite{JY3} and has emerged as a standard condition to derive KKT-type optimality conditions for bilevel optimization problems which are treated by the lower-level value function reformulation \cite{SD2,SD1,SD3,AF2,JY3}.

\begin{definition}
Let $(\bar{x},\bar{y})$ solve \eqref{LLVFR}. Then, \eqref{LLVFR} is said to be \textit{partially calm} at $(\bar{x},\bar{y})$ if there exist constants $r>0$ and $\mu>0$ such that, for all $\alpha \in [-r,r]$ and for all $(x,y) \in B((\bar{x},\bar{y}),r)$ with
\[ 
G(x,y)\leq 0,\quad \hat{g}(x,y)\leq 0,\quad \hat{f}(x,y) -\varphi(x) + \alpha \leq 0, 
\] 
it follows that 
\[ 
F(x,y) - F(\bar{x},\bar{y}) + \mu |\alpha| \geq 0. 
\]
\end{definition}

Note that the definition of partial calmness can also be stated for $(\bar{x},\bar{y})$ being a local minimizer of \eqref{LLVFR}.
Unfortunately, partial calmness is a condition that only holds under restrictive assumptions, for example when the lower-level problem is fully linear and no constraints occur in the upper level \cite[Proposition 4.1]{JY3}. Especially, partial calmness does not necessarily hold as soon as the lower-level problem is linear only with respect to variable~$y$ \cite{PM5}, see Example \ref{ex:partial_calmness} below.
In comparison, on the one hand, LLVFR-subMFC holds for fully linear lower-level problems (possibly with additional concave coupled upper level constraints) as well. Indeed, in this case the lower-level value function is convex \cite[Lemma 2.1]{TT1}, and thus reformulation \eqref{LLVFR} has only concave inequality constraints. Hence, by Theorem \ref{thm:quasinorm} and the fact that concave inequality constraints imply quasinormality \cite[Proposition 3]{JY11}, LLVFR-subMFC automatically holds. On the other hand, we have seen in Section \ref{sec:other_conds} that in situations with even weaker assumptions than concavity of the inequalities, subMFC may hold as well, such that LLVFR-subMFC can be regarded as a promising candidate for an alternative to partial calmness in bilevel optimization. In particular, in the following example from \cite{PM5}, LLVFR-subMFC is satisfied whereas partial calmness does not hold.

\begin{example}[Partial calmness does not imply LLVFR-subMFC]
\label{ex:partial_calmness}
Consider the following instance of an optimistic bilevel optimization problem \eqref{BP}:
\begin{alignat*}{2}
\min_{x,y=(y_1,y_2)} \,& F(x,y)\coloneqq xy_1 \\
\text{ s.t. }& y\in S(x) \coloneqq \argmin_y \{ \hat{f}(x,y)\coloneqq -x^2y_2 \mid\,\, && \hat{g}_1(x,y)\coloneqq y_2 \leq 0,\,\, \\
&&& \hat{g}_2(x,y)\coloneqq -xy_1+y_2 \leq 0 \}. 
\end{alignat*}
As shown in \cite{PM5}, the point $(\bar{x},\bar{y}) \coloneqq (0,(0,0))$ is a global minimizer of the problem and the lower-level value function $\varphi$ is given by
$
\varphi(x) = 0$ for all $x\in\mathbb{R}$. 
Further, it can be verified that this bilevel optimization problem is not partially calm at $(\bar{x},\bar{y})$ \cite{PM5}. In contrast, we now show that LLVFR-subMFC holds at $(\bar{x},\bar{y})$. Concerning part \ref{subMFCII_bilevel} of LLVFR-subMFC, the conditions \eqref{eq:akkt1_bilevel}--\eqref{eq:akkt5_bilevel} for this problem require $(x^k, y^k, \hat{\delta}^k, \gamma_k, \epsilon^k) \to (\bar{x},\bar{y}, 0, 0, 0)$ as $k\to\infty$ and, for all $k\in\mathbb{N}$,
\begin{gather*}
\epsilon^k = \begin{pmatrix} y_1^k \\ x^k \\ 0 \end{pmatrix} + v^k_1  \begin{pmatrix} 0 \\ 0 \\ 1 \end{pmatrix} + v^k_2 \begin{pmatrix} -y_1^k \\ -x^k \\ 1 \end{pmatrix} + w^k  \begin{pmatrix} -2x^ky_2^k \\ 0 \\ -(x^k)^2 \end{pmatrix}, \\
\begin{alignat*}{4}
&& v^k_1 \geq 0,&& y_2^k \leq \hat{\delta}^k_1,&& v^k_1 (y_2^k - \hat{\delta}^k_1) = 0, \phantom{aaaaaaaaaa} \\
&& v^k_2 \geq 0,&&\quad -x^ky_1^k+y_2^k \leq \hat{\delta}^k_2,&&\quad v^k_2 (-x^ky_1^k+y_2^k - \hat{\delta}^k_2 ) = 0, \phantom{aaaaaaaaaa}  \\
\vspace{-2em}
&& w^k \geq 0,&& -(x^k)^2y_2^k \leq \gamma_k,&& w^k (-(x^k)^2y_2^k - \gamma_k ) = 0. \phantom{aaaaaaaaaa}
\end{alignat*}
\vspace{-2em}
\end{gather*}
It can be checked that the sequence 
\[ 
\{ (x^k,y^k,v^k,w^k,\hat{\delta}^k,\gamma_k,\epsilon^k) \}_{k=1}^\infty \coloneqq \left\{ \left( 0,\left(0,-\frac{1}{k}\right),\left(0,0\right),0,(0,0),\frac{1}{k},(0,0,0) \right) \right\}_{k=1}^\infty 
\]
is feasible for both \eqref{LLVFR} and \eqref{LLVFe} and fulfills the above system. As no upper-level constraint exists, $I^{\hat{g}}(x^k,y^k,\hat{\delta}^k)=\emptyset$ and $I^{\hat{f}-\varphi}(x^k,y^k,\gamma_k)=\emptyset$, both conditions \ref{subMFCI_bilevel} and \ref{subMFCII_bilevel} of LLVFR-subMFC are satisfied when setting $I_1\coloneqq\emptyset$, $I_2\coloneqq\emptyset$ and $I_3\coloneqq\emptyset$.
\end{example}

\section{Final remarks}
We introduced subMFC as a problem qualification for Lipschitzian optimization problems which ensures that a feasible point is a KKT point. A comparison with several CQs and the problem qualification calmness revealed that subMFC is a quite weak condition, as it does not imply any of these conditions while being implied by quasinormality. Further, we demonstrated the application of subMFC to bilevel optimization problems and derived the tailored version LLVFR-subMFC that does not imply the well-known partial calmness condition. 
In our future research we aim at applying our findings to different problem classes like MPCCs. Moreover, from a practical standpoint, it can be promising to focus on exploiting our approach when designing numerical algorithms for the treatment of problems belonging to the aforementioned problem classes.

\section*{Acknowledgments}
We would like to thank Patrick Mehlitz for helpful comments on this paper, in particular for providing us with the instance used in Example \ref{ex:strict_implications}.

\section*{Funding}
The work of the first two authors was supported by the Volkswagen Foundation under Grant No. 97 775. The work of the third author was partly funded by the EPSRC grants with references EP/V049038/1 and EP/X040909/1. Moreover, we would like to acknowledge the support of the University of Southampton for a research stay.

\section*{Disclosure of interest}
The authors report there are no competing interests to declare.

\bibliographystyle{spmpsci}    
\bibliography{references}

\vspace{-0.25em}
\appendix
\section{Examples relating subMFC to other conditions}

\subsection{SubMFC and the local error bound condition are independent conditions}
\label{sec:App_subMFC_errorbound}
By the following example, it is first shown that subMFC can hold even if the local error bound condition is violated.

\begin{example}
\label{ex:mehrelementiges_subdiff_error_bound}
Recall the optimization problem from Example \ref{ex:mehrelementiges_subdiff}. The feasible set $D = \{ (x_1,x_2) \in\mathbb{R}^2 \mid g_1(x_1,x_2) \leq 0 \}$ of this problem is visualized in Figure \ref{fig:ex_feasible_region}.

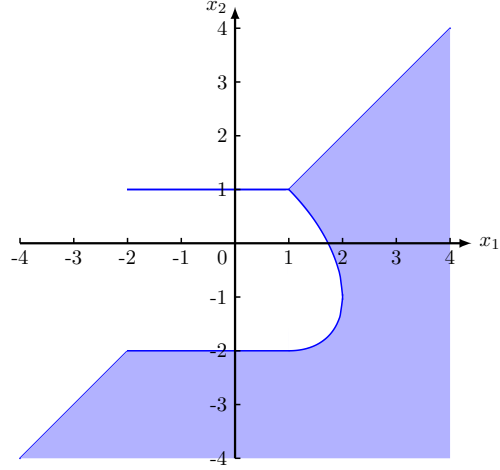
\begin{figure}[h]
	\centering
	\resizebox{0.45\textwidth}{!}{
    \begin{tikzpicture}
    [   axis/.style={black, very thick, -latex, },
        constraint/.style={smooth,blue, thick  },         ]
       
        \draw[domain=-4.01:-1.99, constraint, variable=\x, name path=A] plot ({\x}, {\x});
        \draw[domain=-2.01:1.01, constraint, variable=\x, name path=B] plot ({\x}, 1);
        \draw[domain=-2.01:1.01, constraint, variable=\x, name path=C] plot ({\x}, -2);
        \draw[domain=0.99:4.01, constraint, variable=\x, name path=D] plot ({\x}, {\x});
        \draw[domain=1:2, constraint, variable=\x, name path=E] plot ({\x}, {1/2*(sqrt(3)*sqrt(4-\x^2)-\x)});
        \draw[domain=1:2, constraint, variable=\x, name path=F] plot ({\x}, {1/2*(-sqrt(3)*sqrt(4-\x^2)-\x)});
        	\draw[domain=-4.01:4.01, constraint, variable=\x, name path=G, draw=none] plot ({\x}, -4);
        	\draw[domain=-4.01:4.01, constraint, variable=\x, name path=H, draw=none] plot coordinates {(4,-4) (4,4)};
        	
       	\path [name intersections={of=A and C,by=Q}];
        \path [name intersections={of=C and F,by=R}];
        \path [name intersections={of=F and E,by=S}];
        \path [name intersections={of=D and E,by=T}];
        \path [name intersections={of=D and H,by=U}];
        \path [name intersections={of=G and H,by=V}];
        \path [name intersections={of=G and A,by=W}];
        	
        \colorlet{LightBlue}{white!70!blue}
        	
        	\fill[color=LightBlue]
            plot
            coordinates {
                (W) (V) (U)
            };             
         \fill[color=white]    
            plot[smooth, tension=0]
            coordinates {
                ([xshift= -1pt]Q) ([xshift= 0.1pt]R) ([xshift= -1pt]T)
            };        
		\fill[color=white]    
            plot[smooth, tension=0]
            coordinates {
                (0.99,-1.93) (R) (1.25,-1.977) (1.5,-1.896) (1.75,-1.715) (1.9,-1.495) (1.95,-1.36) (2,-1) (1.95,-0.59) (1.9,-0.41) (1.75,-0.036) (1.5,0.396) (T) (0.99,0.99) (0.95,0.99)
            };  

        \draw[domain=-2.01:1.01, constraint, variable=\x] plot ({\x}, 1);
        \draw[domain=-2.01:1.01, constraint, variable=\x] plot ({\x}, -2);
        \draw[domain=1:2, constraint, variable=\x] plot ({\x}, {1/2*(sqrt(3)*sqrt(4-\x^2)-\x)});
        \draw[domain=1:2, constraint, variable=\x] plot ({\x}, {1/2*(-sqrt(3)*sqrt(4-\x^2)-\x)});
        \draw[domain=2:2, constraint, variable=\x] plot coordinates {(2,-0.95) (2,-1.05)};    
		
		\draw[axis] (-4,0) -- (4.4,0) node[right] {$x_1$}; 
        \draw[axis] (0,-4) -- (0,4.4) node[left] {$x_2$};

		\foreach \x in {-4,-3,-2,-1,1,2,3,4}
    			{   \draw[thick] (\x,0.1) -- (\x,0) node[below] {\x};  	
    				\draw[thick] (0.1,\x) -- (0,\x) node[left] {\x}; 		
    			}
		\draw[thick] (0,0.2) -- (0,0) node[below left] {0};

    \end{tikzpicture}
    }
    \caption{The feasible set $D$ of the optimization problem in Example \ref{ex:mehrelementiges_subdiff_error_bound} is depicted in blue.}
    \label{fig:ex_feasible_region}
\end{figure}
\noindent
We have seen that subMFC holds at the global minimizer $\bar{x}=(\bar{x}_1,\bar{x}_2)=(-2,1)$. However, the local error bound condition does not hold at $\bar{x}$, i.e., we cannot find $c>0$ and $r>0$ such that 
\vspace{-0.7em}
\[ 
d_D(x) \leq c{g_1}_+(x) \quad\forall x=(x_1,x_2)\in B((\bar{x}_1,\bar{x}_2),r). 
\]
To prove this claim, assume that some $r>0$ is fixed. Then, for $x=(x_1,x_2) = (-2, 1+\tau)$ with $\tau\in[0,r]$ we have $x\in B((\bar{x}_1,\bar{x}_2),r)$. 
We obtain $d_D(x) = \tau$, as the point in $D$ that is closest to $x$ is $(-2,1)$, see Figure \ref{fig:ex_feasible_region}. Further, we have 
{\setlength{\belowdisplayskip}{8pt}
\setlength{\abovedisplayskip}{8pt}
\[
{g_1}_+(x) = \max\{0,(1+\tau)^3-3(1+\tau) +2 \} = \tau^3 + 3\tau^2. 
\]
}
Thus, for the local error bound condition to hold, a constant $c>0$ would be needed such that
$
\tau \leq c \,(\tau^3 + 3\tau^2)
$
is satisfied for all $\tau\in[0,r]$.
Obviously, such a constant $c$ does not exist. As $r>0$ was chosen arbitrarily, the local error bound condition cannot hold at $\bar{x}$.
\end{example}

On the other hand, we can also show that the local error bound condition can hold even if subMFC is violated, meaning that both conditions are independent of each other.

\begin{example}
Let us again consider the optimization problem from Example \ref{ex:sub_mfc_kkt}. There, we have seen that subMFC is not valid at the local minimizer $\bar{x}=(\bar{x}_1,\bar{x}_2)=(0,0)$.
In order to show that the local error bound condition holds at $\bar{x}$, we rely on the sufficient condition from \cite[Corollary 2]{AF3}. To apply this condition, we define the function $F(u)\coloneqq g_+(x)$ for which we have $F^{-1}(0)=\bar{x}$ and 
$
\widehat{T}_{F^{-1}(0)}(\bar{x}) = \{ (0,0) \} = \{ d\in\mathbb{R}^2 \mid F'(\bar{x};d)=0\}
$, where $\widehat{T}_{F^{-1}(0)}(\bar{x})$ is a subset of the contingent cone defined in Section \ref{sec:acq_gcq}; see \cite{AF3} for details. Thus, all requirements of \cite[Corollary 2]{AF3} are met and the local error bound condition is satisfied at $\bar{x}$.
\end{example}

\subsection{SubMFC and the cone-continuity property are independent conditions}
\label{sec:App_subMFC_CCP}
The next example demonstrates that CCP can be violated while subMFC holds.

\begin{example}
Consider the following instance of optimization problem \eqref{P}:
\begin{align*}
\min\, (x_1+x_2)x_1 \quad \text{ s.t. } \quad & g_1(x_1,x_2) \coloneqq x_1x_2 \leq 0, \\
& g_2(x_1,x_2) \coloneqq -x_1x_2 \leq 0, \\
& g_3(x_1,x_2) \coloneqq -x_1 \leq 0.
\vspace{-0.5em}
\end{align*}

By the first two constraints, $x_1$ or $x_2$ must be equal to 0 for any feasible point. In the case $x_1=0$, the objective function value becomes $f(0,x_2)=0$, while in the case $x_2=0$, the objective function value is nonnegative, as $f(x_1,0)=x_1^2$. Therefore, all points $(0,x_2)$ with $x_2\in\mathbb{R}$ are global solutions. We will particularly look at the point $\bar{x}=(\bar{x}_1,\bar{x}_2)=(0,1)$.

\begin{enumerate}
\item First, we show that CCP is violated at $\bar{x}$. For this purpose, due to $I(\bar{x})=\{1,2,3\}$, the set $K(\bar{x})$ can be computed as
\begin{align*}
K(\bar{x}) &=  \left\{ \sum_{j \in I(\bar{x})} v_j \nabla g_j(\bar{x}) \mid v_j\geq 0 \right\} \\
&= \left\{ v_1 \begin{pmatrix} \bar{x}_2 \\ \bar{x}_1 \end{pmatrix} + v_2 \begin{pmatrix} -\bar{x}_2 \\ -\bar{x}_1 \end{pmatrix} + v_3 \begin{pmatrix} -1 \\ 0 \end{pmatrix} \mid v_1,v_2,v_3\geq 0 \right\} \\
&= \left\{ \begin{pmatrix} v_1-v_2-v_3 \\ 0 \end{pmatrix} \mid v_1,v_2,v_3\geq 0 \right\} \\
&= \mathbb{R} \times \{ 0 \}.
\end{align*}
For $x^k= (x_1^k,x_2^k)=\left(\frac{1}{k},1\right)$ with $x^k\to \bar{x}$ and $(v_1,v_2,v_3)=(2k,k,k)$ for $k>0$, we get
\[
w^k \coloneqq \sum_{j \in I(\bar{x})} v_j \nabla g_j(x^k) 
= 2k \begin{pmatrix} 1 \\ \frac{1}{k} \end{pmatrix} + k \begin{pmatrix} -1 \\ -\frac{1}{k} \end{pmatrix} + k \begin{pmatrix}-1 \\ 0 \end{pmatrix}= \begin{pmatrix} 0 \\ 1 \end{pmatrix}.
\]
Thus, we obtain $w^k \in K(x^k)$ and $w^k \to \bar{w} \coloneqq (0,1) \in \limsup_{x^k \to \bar{x}} K(x^k)$. 
However, $\bar w\notin K(\bar x)$, such that inclusion \eqref{eq:CCP} does not hold. Hence, CCP is violated at $\bar{x}$.
\item In contrast, subMFC is fulfilled at $\bar{x}$. To see this, consider the AKKT conditions in $\bar{x}$, which require $(x^k, \delta^k,\epsilon^k) \to (\bar{x},0,0)$ as $k\to\infty$ and
\begin{gather*}
\begin{pmatrix}\epsilon^k_1 \\ \epsilon^k_2 \end{pmatrix} = \begin{pmatrix} 2x_1^k+x_2^k \\ x_1^k \end{pmatrix} + u^k_1  \begin{pmatrix} x_2^k \\ x_1^k \end{pmatrix} +
u^k_2  \begin{pmatrix} -x_2^k \\ -x_1^k \end{pmatrix} + u^k_3  \begin{pmatrix} -1  \\ 0 \end{pmatrix}, \\
\begin{alignat*}{4}
&& u^k_1 \geq 0,&& x_1^k x_2^k \leq \delta^k_1,&& u^k_1(x_1^k x_2^k-\delta^k_1)=0, \phantom{aaaaaaaaaaa}\\
&& u^k_2 \geq 0,&&\quad -x_1^k x_2^k \leq \delta^k_2,&&\quad u^k_2(-x_1^k x_2^k-\delta^k_2)=0, \phantom{aaaaaaaaaaa} \\
&& u^k_3 \geq 0,&& -x_1^k \leq \delta^k_3, && u^k_3(-x_1^k-\delta^k_3) = 0 \phantom{,aaaaaaaaaaa}
\end{alignat*}
\end{gather*}
for all $k\in\mathbb{N}$.
One sequence to fulfill the above conditions is given by
\vspace{-0.2em}
\[ 
\{({x}^k,u^k, \delta^k, \epsilon^k)\}_{k=1}^\infty \coloneqq \left\{\left(\left(\frac{1}{k},1\right),\left(0,1,\frac{2}{k}\right),\left(\frac{2}{k},-\frac{1}{k},-\frac{1}{k}\right),\left(0,0\right)\right)\right\}_{k=1}^\infty. 
\]
For this sequence it holds that $I(x^k,\delta^k) = \{ 2,3 \}$ for all $k\in\mathbb{N}$. Hence, with $I\coloneqq \{ 2,3 \}$,\vspace{-0.2em}
\[ 
0 \displaystyle \neq \sum_{j\in {I}} v_j \nabla g_j(\bar{x}) = v_2 \begin{pmatrix} -\bar{x}_2 \\ -\bar{x}_1 \end{pmatrix} + v_3 \begin{pmatrix} -1 \\ 0 \end{pmatrix} = v_2 \begin{pmatrix} -1 \\ 0 \end{pmatrix} + v_3 \begin{pmatrix} -1 \\ 0 \end{pmatrix} 
\]
holds for all $(v_2,v_3) \neq (0,0)$ with $v_2\geq 0,v_3 \geq 0$, i.e., subMFC is satisfied at $\bar{x}$.
\end{enumerate}
\end{example}

Returning to a known example, we can also show the opposite, i.e., that CCP can be valid even if subMFC is violated.

\begin{example}
We again consider Example \ref{ex:sub_mfc_kkt}, where subMFC is violated at the local minimizer $\bar{x}=(\bar{x}_1,\bar{x}_2)=(0,0)$. To verify CCP, we compute $I(\bar{x})=\{1,2,3,4\}$ and
\begin{align*}
K(\bar{x}) &=  \left\{ \sum_{j \in I(\bar{x})} v_j \partial g_j(\bar{x}) \mid v_j\geq 0 \right\} \\
&= \left\{ v_1 \begin{pmatrix} [-1,1] \\ -1 \end{pmatrix} + v_2 \begin{pmatrix} [-1,1] \\ 1 \end{pmatrix} + v_3 \begin{pmatrix} 0 \\ 2 \end{pmatrix}  + v_4 \begin{pmatrix} 0 \\ -2 \end{pmatrix} \mid v_1,v_2,v_3,v_4 \geq 0 \right\} = \mathbb{R}^2.
\end{align*}
Thus, CCP obviously holds at $\bar{x}$.
\end{example}

\subsection{SubMFC and the Guignard constraint qualification are independent conditions}
\label{sec:AppA_subMFC_ACQGCQ}
We now show that subMFC can hold even if both ACQ and GCQ are violated.

\begin{example}
Consider the following instance of optimization problem \eqref{P}:
\begin{align*}
\min \, x_2 \quad \text{ s.t. }\quad & g_1(x_1,x_2) \coloneqq x_1x_2 \leq 0, \\
& g_2(x_1,x_2) \coloneqq -x_1x_2 \leq 0, \\
& g_3(x_1,x_2) \coloneqq x_1^2-x_2 \leq 0.
\end{align*}
The feasible set is given by $D= \{ 0 \} \times [0,\infty)$. Hence, the global solution to the problem is $\bar{x}=(\bar{x}_1,\bar{x}_2)=(0,0)$.

\begin{enumerate}
\item First, we show that both ACQ and GCQ are violated at $\bar{x}$.
The computation of the contingent cone and the cone of locally constrained directions to $D$ at $\bar{x}$ yields
\begin{align*}
T_D(\bar{x}) &= \left\{ \begin{pmatrix} 0 \\ d_2 \end{pmatrix} \in \mathbb{R}^2 \mid d_2 \geq 0 \right\} = D,
\quad L_D(\bar{x}) = \left\{ \begin{pmatrix} d_1 \\ d_2 \end{pmatrix} \in \mathbb{R}^2 \mid d_1\in\mathbb{R},\, d_2 \geq 0 \right\},
\end{align*}
such that ACQ is violated at $\bar{x}$ due to $L_D(\bar{x}) \not\subseteq T_D(\bar{x})$. Due to 
$ \cl \conv T_D(\bar{x}) = T_D(\bar{x}), $
GCQ does not hold at $\bar{x}$ as well.
\item Nevertheless, subMFC holds at $\bar{x}$. To this end, consider the AKKT conditions in $\bar{x}$, which consist of $(x^k, \delta^k,\epsilon^k) \to (\bar{x},0,0)$ as $k\to\infty$ and
\begin{gather*}
\begin{pmatrix}\epsilon^k_1 \\ \epsilon^k_2 \end{pmatrix} = \begin{pmatrix} 0 \\ 1 \end{pmatrix} + u^k_1  \begin{pmatrix} x_2^k \\ x_1^k \end{pmatrix} + u^k_2  \begin{pmatrix} -x_2^k \\ -x_1^k \end{pmatrix} + u^k_3  \begin{pmatrix} 2x_1^k  \\ -1 \end{pmatrix}, \\
\begin{alignat*}{4}
&& u^k_1 \geq 0,&& x_1^k x_2^k \leq \delta^k_1,&& u^k_1(x_1^k x_2^k-\delta^k_1)=0, \phantom{aaaaiiaaa} \\
&& u^k_2 \geq 0,&& -x_1^k x_2^k \leq \delta^k_2,&& u^k_2(-x_1^k x_2^k-\delta^k_2)=0, \phantom{aaaiaiaaa} \\
&& u^k_3 \geq 0,&&\quad (x_1^k)^2-x_2^k \leq \delta^k_3,&&\quad u^k_3( (x_1^k)^2-x_2^k-\delta^k_3) = 0 \phantom{,aaaaaiaia}
\end{alignat*}
\end{gather*}
for all $k\in\mathbb{N}$.
Clearly, the sequence
\[ 
\{(x^k,u^k,\delta^k,\epsilon^k)\}_{k=1}^\infty \coloneqq \left\{(0,0),\left(0,0,1\right),\left(\frac{1}{k},\frac{1}{k},0\right),(0,0)\right\}_{k=1}^\infty
\]
fulfills the above AKKT conditions. For this sequence we have $I(x^k,\delta^k)=\{ 3 \}$ for all~$k\in\mathbb{N}$. Thus, setting $I\coloneqq \{ 3 \}$, it follows that
\vspace{-0.5em}
\[ 
0 \displaystyle \neq \sum_{j\in {I}} v_j \nabla g_j(\bar{x}) = v_3 \begin{pmatrix} 0 \\ -1 \end{pmatrix} 
\]
is fulfilled for all $v_3 > 0$. Hence, subMFC is satisfied at $\bar{x}$.
\end{enumerate}
\end{example}

Again revisiting Example \ref{ex:sub_mfc_kkt}, we see that ACQ and GCQ can be satisfied even if subMFC is not.

\begin{example}
Recall again the optimization problem from Example \ref{ex:sub_mfc_kkt} where we have seen that subMFC does not hold at the local minimizer $\bar{x}=(\bar{x}_1,\bar{x}_2)=(0,0)$.
The computation of the contingent cone and the outer linearization cone to $D$ at $\bar{x}$ yields
\vspace{-0.2em}
\[
T_D(\bar{x}) = \{ 0 \} = L_D(\bar{x}),
\]
such that, due to 
$ \cl \conv T_D(\bar{x}) = T_D(\bar{x}), $
both ACQ and GCQ are fulfilled.
\end{example}

\subsection{SubMFC does not imply calmness in the sense of Clarke}
\label{sec:AppC_subMFC_calmness}
Finally, the following example demonstrates that subMFC can hold even if calmness in the sense of Clarke does not.

\begin{example}
Consider the following instance of optimization problem \eqref{P}:
\vspace{-0.3em}
\begin{align*}
\min \,& f(x_1,x_2) \coloneqq \ln(x_1)-2x_1+x_2 \\
\text{ s.t. }& g_1(x_1,x_2) \coloneqq 2x_1- x_2^3+9 x_2^2- 29x_2+31 \leq 0, \\
& g_2(x_1,x_2) \coloneqq -\frac{1}{2}x_1+\frac{1}{54} x_2^3- \frac{1}{6} x_2^2+ x_2-\frac{3}{2} \leq 0, \\
& g_3(x_1,x_2) \coloneqq x_2-3 \leq 0.
\end{align*}

\begin{enumerate}
\item We first show that $\bar{x}=(\bar{x}_1,\bar{x}_2)=(1,3)$ is a global solution to the problem. For this purpose, note that $g_1$ and $g_2$ provide an upper and a lower bound for $x_1$. Indeed, for any feasible point $x=(x_1,x_2)$, it is required that
\vspace{-0.5em}
\[ 
\frac{1}{27} x_2^3- \frac{1}{3} x_2^2+ 2x_2-3 \leq (x_1 \leq)\, \frac{1}{2} x_2^3-\frac{9}{2} x_2^2+ \frac{29}{2}x_2-\frac{31}{2} 
\]
holds, which is equivalent to 
\vspace{-0.3em}
\[
\frac{25}{54} \left(x_2^3-9x_2^2+27x_2-27\right) \geq 0 \Longleftrightarrow \frac{25}{54} \left(x_2-3\right)^3 \geq 0 \Longleftrightarrow x_2-3 \geq 0.
\]
Together with $g_3(x_1,x_2)\leq 0$, this yields $x_2=3$ for any feasible point. In turn, $g_1$ and $g_2$ for $x_2=3$ read as
\vspace{-0.3em}
\[ 
g_1(x_1,3) = 2x_1-2 \,\,\,\,\, \text{and} \,\,\,\,\, g_2(x_1,3) = -\frac{1}{2}x_1 + \frac{1}{2}. \vspace{-0.5em}
\]
Thus, it can be verified that the point $\bar{x}= (1,3)$ is the only feasible point and therefore a global solution to the problem with $f(1,3) = 1$.

\item We now demonstrate that calmness is violated at $\bar{x}=(1,3)$.
To this end, let us suppose the contrary, i.e., that calmness holds at the solution $\bar{x}= (1,3)$. Then, there exist $r>0$ and $\mu>0$ such that for all $\alpha \in B(0,r)$ it holds that
\begin{equation} 
\ln(x_1)-2x_1+x_2- 1 + \mu \norm{\alpha} \geq 0 \quad \forall x \in \{ (x_1,x_2) \in B(\bar{x},r) \mid g(x_1,x_2) + \alpha \leq 0\}. \label{eq:calm_ex}
\end{equation}
Consider the sequence $\{ (\alpha^k, x^k) \}_{k=1}^\infty$ with
\[ 
(\alpha^k, x^k) \coloneqq \left( \left(-\frac{1}{k^3},0,0\right), \left( 1-\frac{1}{k}, 3-\frac{1}{k} \right) \right)
\]
which fulfills $(\alpha^k, x^k) \to (\left(0, 0,0\right), (\bar{x}_1,\bar{x}_2))$ as $k\to\infty$ and
\[
g_1(x_1^k,x_2^k) + \alpha^k_1 = 0,\,\, g_2(x_1^k,x_2^k) + \alpha^k_2 = -\frac{1}{54k^3}<0,\,\, g_3(x_1^k,x_2^k) + \alpha^k_3 = -\frac{1}{k} < 0.
\]
However, as $\lVert \alpha^k \rVert =  \frac{1}{k^3}$, it follows from \eqref{eq:calm_ex} that
\vspace{-0.5em}
\[
 \ln(x_1^k)-2x_1^k+x_2^k-1 + \mu \lVert \alpha^k \rVert = \frac{1}{k} + \ln\left(1-\frac{1}{k}\right) + \frac{\mu}{k^3} \geq 0. \vspace{-0.4em}  
\]
This inequality is equivalent to
\vspace{-0.7em}
\begin{equation}
\mu\ge -k^2\left(k\ln\left(1-\frac1k\right)+1\right), \label{eq:calm_ex2}
\end{equation}
and it can be shown that the right-hand side of \eqref{eq:calm_ex2} converges to $\infty$ for $k\to\infty$.
Thus, we cannot find a finite $\mu>0$ to fulfill the inequality, and hence, the problem is not calm at $\bar{x}$.

\item In contrast, subMFC is fulfilled at $\bar{x}$. To verify this claim, we take a look at the AKKT conditions in $\bar{x}$, i.e., we need to satisfy $(x^k, \delta^k,\epsilon^k) \to (\bar{x},0,0)$ as $k\to\infty$ and
\vspace{-0.3em}
\begin{gather*}
\hspace{-0.5em}
\begin{pmatrix}\epsilon^k_1 \\ \epsilon^k_2\end{pmatrix}
= \begin{pmatrix}\frac{1}{x_1^k} -2 \\ 1 \end{pmatrix} + u^k_1  \begin{pmatrix} 2 \\ -3(x_2^k)^2+18x_2^k-29 \end{pmatrix} + u^k_2 \begin{pmatrix} -\frac{1}{2} \\ \frac{1}{18}(x_2^k)^2 -\frac{1}{3}x_2^k+1 \end{pmatrix} + u^k_3 \begin{pmatrix} 0 \\ 1 \end{pmatrix},\\
\begin{alignat*}{4}
&& u^k_i \geq 0,&&\quad g_i(x_1^k,x_2^k) \leq \delta^k_i,&&\quad u^k_i (g_i(x_1^k,x_2^k) - \delta^k_i) = 0, \quad\forall i=1,\dotsc,3 \phantom{aaaaa}
\end{alignat*}
\end{gather*}
for all $k\in\mathbb{N}$.
Consider the sequence
\[ 
\hspace{-0.4em}
\left\{ ({x}^k,u^k, \delta^k, \epsilon^k) \right\}_{k=1}^\infty \coloneqq 
 \left\{ \left(\left(1-\frac{1}{k},3-\frac{1}{k}\right), 
 \left(\frac{k-2}{2k-2},0,0\right), 
 \left(\frac{1}{k^3},0,0\right), 
 \left(0,\epsilon^k_2\right)\right) \right\}_{k=1}^\infty
\]
with $\epsilon^k_2 \coloneqq \frac{2k^2-3k+6}{2(k-1)k^2} \to 0$ as $k \to\infty$.
Then, the first line of the AKKT system reads as
\[ 
\begin{pmatrix} 0 \\ \epsilon^k_2 \end{pmatrix} = \begin{pmatrix} \frac{k}{k-1} -2 \\ 1 \end{pmatrix} + 
\frac{k-2}{2k-2} \cdot \begin{pmatrix} 2 \\ -\frac{3}{k^2}-2 \end{pmatrix},
\]
while the other AKKT conditions are fulfilled as well due to $(x^k, \delta^k,\epsilon^k) \to (\bar{x},0,0)$ as $k\to\infty$ and
\[
g_1(x_1^k,x_2^k) = \frac{1}{k^3} = \delta^k_1,\quad g_2(x_1^k,x_2^k) =-\frac{1}{54k^3} < 0 = \delta^k_2,\quad g_3(x_1^k,x_2^k)= -\frac{1}{k} < 0 = \delta^k_3.
\]
Thus, for the sequence $\{(x^k,u^k,\delta^k,\epsilon^k)\}_{k=1}^{\infty}$, it holds for all $k\in\mathbb{N}$ that $I(x^k,\delta^k) = \{ 1 \}$.
Further, with $I\coloneqq \{ 1 \}$, condition \ref{subMFCI} of subMFC means that
\[ 
0 \displaystyle \neq \sum_{j\in {I}} v_j \nabla g_j(\bar{x}) = v_1  \begin{pmatrix} 2 \\ -3\bar{x}_2^2+18\bar{x}_2-29 \end{pmatrix} = v_1  \begin{pmatrix} 2 \\ -2 \end{pmatrix}  
\]
has to be fulfilled for all $v_1 > 0$, and this is indeed the case. Thus, subMFC holds at~$\bar{x}$.
\end{enumerate}
\end{example}

\section{Examples and results verifying relations between different conditions}
\label{sec:examples_strictness}
In this section, we briefly summarize how the relations depicted in Figure \ref{fig:CQ_implications}, which have not yet been discussed in the course of this paper, can be verified. Results concerning implications can be obtained as follows: 
\begin{itemize}
\item The implication between LICQ and MFCQ follows trivially from the respective definitions, as LICQ at a feasible point $\bar{x}$ of \eqref{P} demands the linear independence of the subgradients of all constraints that are active in $\bar{x}$, whereas for MFCQ only their positive linear independence is required. Analogously, the implications between LICQ and SMFC as well as SMFC and MFCQ are obtained similarly to the smooth case.
\item The nonsmooth extension of CRCQ is given in \cite{MX1} and trivially yields that CRCQ is implied by LICQ. 
\item In \cite{RJ1}, it is shown that CRCQ and MFCQ do not imply each other.
\end{itemize}
Remaining results concerning the strictness of the presented implications can be established using the following examples:
\begin{itemize}
\item The feasible set given by $g_1(x) \coloneqq 2x \leq 0$ and $g_2(x) \coloneqq x \leq 0$ and the feasible point $\bar{x}=0$ can be used to see that the implications between LICQ and MFCQ as well as LICQ and CRCQ are strict.
\item The optimization problem \eqref{P} with $f(x)\coloneqq x^2$, $g_1(x)\coloneqq x$ and $g_2(x) \coloneqq -x$ with the global optimizer $\bar{x}=0$ can be used to see that the implications between MFCQ and pseudonormality as well as MFCQ and calmness are strict. 
\item Further, \cite[Counterexample 4.4]{RA1} shows that pseudonormality can hold even if not all constraint functions are concave and that the implication between CPLD and quasinormality is strict.
\item For the feasible set given by $g_1(x_1,x_2) \coloneqq -x_2e^{x_1}\leq 0$ and $g_2(x_1,x_2) \coloneqq x_2\leq 0$, which is inspired by \cite{RA1}, the local error bound condition is satisfied at $\bar{x}=(0,0)$ while CPLD is not.
\item When considering the feasible set given by $g_1(x_1,x_2)\coloneqq |x_1|-x_2\leq 0$ and $g_2(x_1,x_2) \coloneqq x_2 \leq 0$ and the point $\bar{x}=(0,0)$, it can be shown that the implications between quasinormality and the local error bound condition as well as between quasinormality and ACQ are strict.
\item In \cite{RA5}, it is shown that the implication between CCP and ACQ is strict.
\item The strictness of the implication between ACQ and GCQ can be verified by considering the feasible set given by $g_1(x_1,x_2)\coloneqq -x_1\leq 0$, $g_2(x_1,x_2)\coloneqq -x_2\leq 0$ and $g_3(x_1,x_2)\coloneqq x_1x_2\leq 0$ and the feasible point $\bar{x}=(0,0)$.
\end{itemize}

\end{document}